\documentclass[centertags,reqno]{amsart}
\usepackage{amsmath,amssymb,amscd}
\newcommand{\N}{{\mathbf N}}
\newcommand{\A}{{\mathbf A}}
\newcommand{\Z}{{\mathbf Z}}

\newcommand{\C}{{\mathbf C}}
\newcommand{\Q}{{\mathbf Q}}
\renewcommand{\P}{{\mathbf P}}

\newcommand{\K}{{\mathbf K}}

\newcommand{\Gam}{\Gamma}
\newcommand{\Lam}{\Lambda}
\newcommand{\lam}{\lambda}
\newcommand{\bet}{\beta}
\newcommand{\del}{\delta}
\newcommand{\sig}{\sigma}
\newcommand{\alp}{\alpha}
\newcommand{\om}{\overline{m}}
\newcommand{\oz}{\overline{z}}
\newcommand{\mt}{\mbox{$\mathrm{mt}$}}
\newcommand{\lra}{\longrightarrow}

\newcommand{\kc}{{\mathcal C}}

\newcommand{\kg}{{\mathcal G}}

\newcommand{\kj}{{\mathcal J}}

\newcommand{\ko}{{\mathcal O}}

\newcommand{\ks}{{\mathcal S}}
\newcommand{\kt}{{\mathcal T}}

\newcommand{\kx}{{\mathcal X}}

\newtheorem{lemma}[equation]{Lemma}
\newtheorem{proposition}[equation]{Proposition}
\newtheorem{theorem}{Theorem}

\theoremstyle{definition}
\newtheorem{definition}[equation]{Definition}

\theoremstyle{remark}
\newtheorem{remark}[equation]{Remark}

\newtheorem{examples}{Examples} 

\begin{document}

\title[Plane curves with prescribed singularities]{Plane
  curves of minimal degree with prescribed singularities}
\author[G.-M. Greuel]{Gert-Martin Greuel}
\address{
Universit\"at Kaiserslautern\\
Fachbereich Mathematik\\
Erwin-Schr\"odinger-Stra\ss e\\
D -- 67663 Kaiserslautern}
\email{greuel@mathematik.uni-kl.de}
\author[C. Lossen]{Christoph Lossen}
\address{
Universit\"at Kaiserslautern\\
Fachbereich Mathematik\\
Erwin-Schr\"odinger-Stra\ss e\\
D -- 67663 Kaiserslautern}
\email{lossen@mathematik.uni-kl.de}
\author[E. Shustin]{Eugenii Shustin}
\address{
Tel Aviv University\\
School of Mathematical Sciences\\
Ramat Aviv\\
ISR -- Tel Aviv 69978}
\email{shustin@math.tau.ac.il}

\begin{abstract}
We prove that there exists  a positive $\alp$ such that
for any integer \mbox{$d\ge 3$} and any topological types
\mbox{$S_1,\dots,S_n$} of plane curve singularities, satisfying
\mbox{$\mu(S_1)+\dots+\mu(S_n)\le\alp d^2$},
there exists a reduced irreducible plane curve of degree $d$
with exactly $n$ singular points of types \mbox{$S_1,\dots,S_n$},
respectively. This estimate is optimal with respect to the
exponent of $d$. In particular, we prove that for any topological type $S$
there exists an irreducible polynomial of degree \mbox{$d\le 14\sqrt{\mu(S)}$}
having a singular point of type $S$.
\end{abstract}

\maketitle

\section*{Introduction}

Throughout the article we consider all objects to be defined over
an algebraically closed field $\K$ of characteristic zero.

In the paper we deal with the following classical problem:
given an integer \mbox{$d\ge 3$} and types \mbox{$S_1,\dots,S_n$} of plane
curve singularities, does there exist a reduced irreducible
plane curve of degree $d$ with exactly $n$ singular points of types
\mbox{$S_1,\dots,S_n$}, respectively? The complete answer is known for
nodal curves \cite{Sev}: an irreducible curve of degree $d$, with $n$
nodes as its only singularities, exists if and only if
$$0\le n\le \frac{(d-1)(d-2)}{2}\: .$$
For other singularities, even for ordinary cusps there is no
complete answer. Namely, various restrictions are found
(from Pl\"ucker formulae to inequalities by
Varchenko \cite{Va} and Ivinskis \cite{HF,Ivi}), which read
\begin{equation} 
  \label{0.1}
\sum_{i=1}^n \sig(S_i) < \alp_2 d^2+\alp_1 d+\alp_0, \;\;\;\alp_2=
\mbox{const}>0, 
\end{equation}
with some positive invariants $\sig$ of singular points which are at most
quadratic in d.
We want to give an {\bf asymptotically optimal} sufficient existence condition,
that is a condition of type
$$\sum_{i=1}^n\sig(S_i)<\alp d^2+o(d^2),\quad \alp\leq\alp_2,$$
providing (\ref{0.1}) is necessary.

Note that an asymptotically exact condition, that is \mbox{$\alp = \alp_2$},
is hardly attainable. For example, there exist
curves of degree \mbox{$d=2\cdot 3^k$}, \mbox{$(k=1,2,\dots)$} with
$9(9^k-1)/8=9d^2/32+O(d)$ ordinary cusps
\cite{H}. But here 
the number of conditions imposed by the cusps is \mbox{$d^2/16+O(d)$} more than
the dimension of the space of curves of degree $d$, therefore one cannot
expect that all intermediate quantities of cusps may be realized.

The only previously known general sufficient condition for the
existence of a curve with given singularities was (see \cite{Sh})
\begin{equation}
  \label{no}
  \sum_{i=1}^n(\mu(S_i)+4)(\mu(S_i)+5)\le\frac{(d+3)^2}{2}\ .
\end{equation}
It is not asymptotically optimal, because
the left--hand side may be about $d^4$.

The goal of this paper is:

\begin{theorem}
  \label{0.2}
For any integer \mbox{$d\ge 1$} and topological types \mbox{$S_1,\dots,S_n$} of
plane curve singularities, satisfying
\begin{equation}
  \label{0.3}
  \sum_{i=1}^n \mu(S_i) \le \frac{d^2}{392}\: ,
\end{equation}
there exists a reduced irreducible plane projective curve of degree
$d$ with exactly $n$ singular points of types \mbox{$S_1,\dots,S_n$},
respectively.
\end{theorem}

This estimate is asymptotically optimal, because always
$$\sum_{i=1}^n \mu(S_i)<d^2\ .$$
The constant in (\ref{0.3}) is not the best possible. Our method could give
a bigger constant, providing more tedious computations. For certain
classes of singularities such as simple or ordinary, there are
much better results (see, for instance, \cite{GLS}, section 3.3 and \cite{Sh}).

The problem is of interest even for one individual singularity.
Given a singularity $S$, what is the minimal degree of a reduced irreducible
plane projective curve having 
this singularity at the origin? The classical upper bound
is the determinacy bound \mbox{$\mu(S)+2$} \cite{To}, whereas a lower bound
is \mbox{$\sqrt{\mu(S)}+1$} (coming from intersecting two generic polars and
B\'ezout's Theorem). We claim

\begin{theorem}
  \label{0.4}
For any topological type $S$ of plane curve singularities there
exists a reduced irreducible plane projective curve of degree \mbox{$\le
  14\sqrt{\mu(S)}$} with singularity $S$ at the origin.
\end{theorem}

We should like to thank the Deutsche Forschungsgemeinschaft and the grant G
039-304.01/95 of the German--Israeli Foundation for financial
support. 

\section{Strategy of the Proof}
\setcounter{equation}{0}

To emphasize what is new in our approach we describe shortly
the main previously known constructions.

The first one is to construct, somehow, a curve of the given degree,
which is degenerate
with respect to the required curve, and then to deform it in order
to obtain the prescribed singularities. For example, Severi \cite{Sev}
showed that singular points of a nodal curve, irreducible or not,
can be smoothed out or preserved independently. Hence, taking the
union of generic straight lines and smoothing out suitable intersection
points, one obtains irreducible curves with any prescribed number of
nodes, allowed by Pl\"ucker's formulae. Attempts to extend this construction
on other singularities give curves with a number of singularities
bounded from above by a linear function of the degree $d$
(see, for example, \cite{GrM} for curves with nodes,
cusps and ordinary triple points), because of the very restrictive
requirement of the independence of deformations of singular points.

The second way consists of a construction especially adapted to
the given degree and given collection of singularities. It may be
based on a sequence of rational transformations of the plane
applied to a more or less simple initial curve in order to get
the required curve.
Or it may
consist in an invention of a polynomial defining the required curve.
This can be illustrated by constructions of singular curves
of small degrees as, for instance, in \cite{Wl1}, \cite{Wl2},
or by the construction of cuspidal curves as in \cite{H}, cited in the
Introduction. Two main difficulties do not allow the appliance of this
approach to a wide class of degrees and singularities: (1) for any
new degree or singularity one has to invent a new construction,
(2) even if one has constructed a curve with a lot of singularities, like in
\cite{H}, it is hard to check that these singular points can be smoothed out
independently and any intermediate numbers of singularities can be realized.

Another idea, based on a modification of the Viro method of gluing
polynomials (see the original method in \cite{Vi})
and on the independence of singular point deformations,
was suggested in \cite{Sh}. This method, from the very beginning, requires
a collection of ``base curves'' with given singularities
(as, for instance, in Theorem \ref{0.4}), which originally provides
only non--optimal results such as (\ref{no}) for arbitrary singularities.

In our proof we use the previous constructions and introduce 
the following new element. With reduced germs of plane curves we associate
a class ${\kg \ks}$ of irreducible zero--dimensional schemes, called below
{\bf generalized singularity schemes}. Further we proceed in three main
steps.

{\it Step 1}. Given a topological type $S$ of plane curve singularities,
we show that there exists a scheme \mbox{$X\in{\kg\ks}$} with \mbox{$\deg
  X\le a_1\mu(S)$}, 
\mbox{$a_1=\mbox{const}>0$} such that the relation
$$h^1(\P^2,\kj_{X/\P^2}(d))=0,$$
where \mbox{$\kj_{X/\P^2}\subset{\ko}_{\P^2}$} is the ideal sheaf of $X$,
suffices for
the existence of a curve of degree $d$ with a singular point of type $S$
(see Lemmas \ref{4.1}, \ref{4.11} below).

{\it Step 2}. For our purposes we have to provide the previous 
$h^1$--vanishing as \mbox{$d\le a_2\sqrt{\deg X}$},
\mbox{$a_2=\mbox{const}>0$}. To do this,
we observe that in the first step, $X$ can be replaced by a generic
scheme $X'$ in the same Hilbert scheme. Then we follow basically
Hirschowitz \cite{HA}, who obtained, in an analogous manner,
the $h^1$--vanishing
for schemes of generic fat points in the plane. Namely, we fix a
straight line $L$ and apply an inductive procedure described
in Sections 3 and 4, which consists
of a passage from $X$ and $d$ to the residue scheme $X:L$ (called below the
{\bf reduction of} $X$) and $d-1$. Each time we have to verify that
$X:L$ belongs to $\kg\ks$ (Proposition \ref{1.11}), and that
$$a_3d\le\deg(X\cap L)\le d+1,\quad a_3=\mbox{const}>0.$$
The latter relation is achieved by means of two operations: {\bf
  specialization}
of the scheme $X$ with respect to $L$ (Lemma \ref{1.14}), and 
{\bf extension} of $X$ (Definitions \ref{1.19}, \ref{1.21}, Lemma \ref{1.20})
when the specialization fails.
 
{\it Step 3}. The final stage is a construction of curves with many
singular points, done by means of a version of the Viro method
(Section 6). Given topological
singularities \mbox{$S_1,\dots,S_n$}, we find curves of degrees
$$d_i\le a_4\sqrt{\mu(S_i)},\quad i=1,\dots ,n,\quad a_4=\mbox{const}>0,$$
each having a singular point of the corresponding type. Then 
we take a curve of degree 
$$d\le a_5\sqrt{d_1^2+...+d_n^2},\quad a_5=\mbox{const}>0,$$
with $n$ generic points of multiplicities \mbox{$d_1,\dots ,d_n$},
respectively, and deform these points in order to obtain
the given singularities on a curve of any degree
$$d\ge a_6\sqrt{\mu(S_1)+...+\mu(S_n)},\quad a_6=\mbox{const}>0.$$

\section{Singularity Schemes, Reductions and Extensions}
\setcounter{equation}{0}

Throughout this section, $S$ denotes a smooth surface, \mbox{$z \in S$}, and
$C$ 
a reduced curve on $S$.   Since our statements are local, we may assume that
$C$, or the germ $(C,z)$, is given by a power series which, by abuse of
notation,  is also
denoted by $C$ or by $(C,z)$.  If \mbox{$z \not\in C$}, then $(C,z)$ denotes
the empty germ or a unit of $\ko_{S,z}$.   Later, $z$ denotes also a finite set
of points of $S$ and $(C,z)$ the corresponding multigerm.

\begin{definition}\label{1.1}
  The {\bf multiplicity} of $C$ at $z$ is the non--negative integer
\[
\mt \,(C,z) = \operatorname{max} \{n \in \Z \mid C \in \mathfrak{m}^n_z\},
\]
where $\mathfrak{m}_z$ is the maximal ideal of $\ko_{S,z}$, the analytic local
ring of $S$ at $z$.

If \mbox{$z \in C$}, we define, as usual, (cf.\ \cite{Z}, \cite{W}, \cite{Te},
\cite{BK}) the {\bf topological type} (or {\bf equisingularity type}) of the
germ $(C,z)$ by the following discrete characteristic:  the embedded
resolution tree of $(C,z)$ and the multiplicities of the total transforms of
$(C,z)$ at infinitely near points (including $z$).
Two germs with the same topological types are called equivalent (notation
$\sim$).
\end{definition}

\begin{definition}\label{1.2}
  $z$ is called an {\bf essential} point of $C$ if \mbox{$z \in C$}, and if the
  germ $(C,z)$ is not smooth.   If \mbox{$z \in C$} and if \mbox{$q \not= z$}
  is infinitely near to 
  $z$, we denote by $C_{(q)}$, respectively $\widehat{C}_{(q)}$, the
  corresponding strict, respectively total, transforms under the composition of
  blowing--ups  \mbox{$\pi_{(q)} : S_{(q)} \lra S$} defining $q$.   We call $q$
  {\bf essential} if it is not 
  a node (ordinary double point) of the union of $C_{(q)}$ with the reduced
  exceptional divisor. 
\end{definition}

We shall introduce now the singularity scheme, respectively the generalized
singularity scheme, of $(C,z)$, which are zero--dimensional subschemes of $S$
and which encode to a certain extent the topological type of $(C,z)$,
respectively together with some higher order tangencies.

For \mbox{$z \in C$} let \mbox{$T(C,z)$} denote the (infinite) complete
embedded resolution tree of $(C,z)$ with vertices the points infinitely near to
$z$. 
It is naturally oriented, inducing a partial ordering on its vertices such
that \mbox{$z < q$} for all \mbox{$q \in T(C,z) \backslash \{z\}$}.   If
\mbox{$z \not\in
  C$} we define $T(C,z)$ to be the empty tree.   Moreover, let
\[
T^\ast (C,z) := \{q \in T(C,z) \mid q \mbox{ is essential}\}
\]
denote the tree of essential points of $(C,z)$, which is a finite subtree of $T
(C,z)$.

\begin{definition}\label{1.3}
  Let \mbox{$T^\ast \!\subset T(C,z)$} be a finite, connected subtree
  containing the essential tree
  $T^\ast(C,z)$.   For any point \mbox{$q \in T^\ast$} and any \mbox{$f \in
  \ko_{S, z}$}
  denote by 
  $f_{(q)}$, respectively $\hat{f}_{(q)}$, the strict, respectively total,
  transform under the modification $\pi_{(q)}$ defining $q$.   Put \mbox{$m_q
  := \mt \,(C_{(q)}, q)$},
  \mbox{$\hat{m}_q := \mt \,(\widehat{C}_{(q)}, q)$} and define the ideal
\[
J := J(C,T^\ast) := \{f \in \ko_{S,z} \mid \mt \,(\hat{f}_{(q)}, q) \ge
\hat{m}_q,\;\; q \in T^\ast\} \subset \ko_{S,z}
\]
and the subscheme of $S$ defined by $J$,
\[
X := X(C,T^\ast) = Z(J),\quad \ko_{X,z} := \ko_{S,z}/J,
\]
which is concentrated on $\{z\}$. $X$ is called a {\bf generalized singularity
  (scheme)} and the class of zero--dimensional subschemes of $S$, constructed
in this way, is denoted by $\kg \ks$.  The subclass of schemes \mbox{$X \in \kg
  \ks$} 
with \mbox{$T^\ast \!= T^\ast(C,z)$} is denoted by $\ks$, \mbox{$X \in \ks$}
  is called a {\bf singularity (scheme)}.
\end{definition}

\begin{examples} 
\begin{enumerate}
\item Let $(C,z)$ be smooth.  If \mbox{$T^\ast \!= \emptyset$}, we obtain
  \mbox{$J = \ko_{S,z}$} and \mbox{$X = \emptyset$}.  If \mbox{$T^\ast \!=
  \{ z \!=\! q_0,q_1,\dots, q_n\}$} and \mbox{$C =y$} with respect to local
  coordinates 
  $(x,y)$ at $z$, then \mbox{$J = \langle y, x^{n+1}\rangle$}.

\item If $(C,z)$ is an ordinary $r$--fold singularity ($r$ smooth branches
  with different tangents) and if \mbox{$T^\ast \!=
  \!\{z\}\;  (= T^\ast(C,z))$} then \mbox{$J = \mathfrak{m}^r_z$}.
\end{enumerate}
\end{examples}

The following lemma shows the relation of $X$ to equisingular deformations of
$(C,z)$.   Note that the germ $(C,z)$ defining $X$ is not uniquely determined
by $X$, and that the tree $T^\ast$ is part of the data of $X$.   We write
$J_X$ and $T^\ast_X$ for the ideal $J$ and the tree $T^\ast$ belonging to $X$.

By the following lemma, though $(C,z)$ is not uniquely determined by $X$, all
topological invariants of $(C,z)$ can be associated uniquely to $X$.   In
particular, we
define the {\bf multiplicity} $\mt \,X$, the {\bf Milnor number} $\mu(X)$ and
the {$\mathbf \delta$}--{\bf invariant} $\delta(X)$ as those of $(C,z)$.

\pagebreak[3]
\begin{lemma}\label{1.4}
  Let \mbox{$X \in \kg \ks$} be a generalized singularity scheme.

  \begin{enumerate}
  \item[(i)] If $X$ is defined by $(C,z)$ then a generic element of $J_X$ is
    topologically equivalent to $(C,z)$.
  \item[(ii)] The set of base points of the ideal $J_X$ is equal to $T^\ast_X$,
    that is, the strict transforms of two generic elements in $J_X$ intersect
    exactly in $T^\ast_X$.
  \end{enumerate}
\end{lemma}

\begin{proof}
Adding a generic element \mbox{$f\in J_X$} to the equation of
\mbox{$(C,z)\subset (S,z)$} defines 
another generic element of $J_X$ having exactly the given multiplicities
$\hat{m}_q$, \mbox{$q\in T^{\ast}$}, in particular, it is  topologically
equivalent to $(C,z)$, since $T^\ast(C,z) \subset T^\ast$. \\
Moreover, this shows that the strict transforms of two generic elements
have the same multiplicity \mbox{$m_q\geq 1$} at \mbox{$q\in T^\ast_X$}. On the
other hand, let  \mbox{$q\in T(C,z)\backslash T^\ast_X$}, $Q$ be the
corresponding branch of $(C,z)$ and $\bar{q}$ the predecessor of $q$
in \mbox{$T(C,z)$}. Then, slightly changing the tangent direction of the strict
transform of $Q$ in $\bar{q}$, blowing down and composing with the other
branches of $(C,z)$ defines an element of $J_X$ whose strict transform at $q$
does not contain $q$.
\end{proof}

The concepts developed so far generalize immediately to multigerms $(C,z)$,
\mbox{$z= \{z_1, \dots, z_k\} \subset S$}. Then $T(C,z)$ and $T^\ast(C,z)$ are
finite unions of trees. 
For \mbox{$T^\ast(C,z) \subset T^\ast \subset T(C,z)$} such that
\mbox{$T^\ast \!\cap T(C,z_i)$} is a finite and connected subtree, we can
define \mbox{$J(C,T^\ast) \subset \ko_{S,z} =
  \prod\nolimits_{i=1}^k 
\ko_{S,z_i}$} and \mbox{$X = X(C,T^\ast) = Z(J)$}, as before.   
$X$ is then a reducible subscheme of $S$, concentrated on \mbox{$z_1, \dots,
  z_k$}.  Let
$\mt \,X,\; \mu(X),\; \delta(X)$ denote the sum of
the corresponding invariants at \mbox{$z_1, \dots, z_k$}.   We need this
generalization after blowing up.

Let \mbox{$z \in C \subset S$} be a point and \mbox{$\widehat{S} \lra S$} be
the blowing--up of $z$.   We denote by $\widehat{C}$, respectively $C^\ast$,
the total, 
respectively strict, transform of $C$, $E$ the reduced exceptional divisor and
\mbox{$\hat{z} := E \cap C^\ast$}.   \mbox{$(C^\ast\!, \hat{z}) \subset
  (\widehat{S}, \hat{z})$} is a (multi)germ.

For any \mbox{$f \in \ko_{S,z}$} satisfying \mbox{$\mt \,(f,z) \ge m := \mt
  \,(C,z)$} we may divide
the total transform $\hat{f}$ by the $m$'th power of $E$ and we shall denote
this multigerm at $\hat{z}$ by \mbox{$\hat{f} : mE$}.   If \mbox{$m =
  \mt \,(f,z)$}, then 
\mbox{$\hat{f} :m E = f^\ast$} is the strict transform of $f$.   Note that for
\mbox{$q \in T(C,z)\backslash \{z\}$}
\[ 
\mt \,(\hat{f}_{(q)}, q) = \mt \,((\hat{f} : mE)_{(q)}, q) + k(q)\cdot m,
\]
where \mbox{$k(q)\in \N$} is independent of $f$. This holds
  especially for 
  \mbox{$f = C$}, hence we obtain for \mbox{$T^\ast(C,z) \subset
  T^\ast \subset T(C,z)$}:

\begin{lemma}\label{1.5}
  With the above notations:
\[
f \in J(C,T^\ast) \Leftrightarrow 
  \left( \mt \,(f,z) \ge m \mbox{ and }
\hat{f} : mE \in J(C^\ast, T^\ast \backslash\{z\})\right).
\]
\end{lemma}

Let us denote by \mbox{$\mt \,(X,q) := m_q$}, \mbox{$q \in T^\ast_X$}, the {\bf
  multiplicity of $\bf X$ at $\bf q$} and by \mbox{deg$(X) :=
  \dim_K(\ko_{S,z}/J_X)$} the {\bf degree of $\bf X$}.

\begin{lemma}\label{1.6}
  For \mbox{$X \in \kg \ks$}, \mbox{$T^\ast \!= T^\ast_X$} and \mbox{$m_q =
  \mt \,(X,q)$} we have 
\[
\deg(X)  \,=\, \sum_{q \in T^\ast}
\frac{m_q(m_q + 1)}{2} \,=\, \delta(X) + \sum_{q \in T^\ast} m_q.
\]
\end{lemma}

\begin{proof}
The second equality follows from 
$$ \delta(C,z) = \sum_{q\in T^\ast} \frac{m_q(m_q-1)}{2}.$$
For the proof of the first one, cf.\ \cite{Casas-Alvero}, Proposition 6.1.
\end{proof}

\begin{lemma}\label{1.7}
  Let \mbox{$X \in \ks$} be defined by a germ $(C,z)$ and let \mbox{$I^{es}
  \subset \ko_{S,z}$} be the equisingularity ideal of $(C,z)$ in the sense of
  Wahl $($\cite{W}, cf.\ also \cite{DH}$)$, then \mbox{$I^{es} \supset J_X$}.
\end{lemma}

\begin{proof}
Clear from the definitions.
\end{proof}

\begin{lemma}\label{1.8}
  Let \mbox{$Q_1, \dots, Q_r$} denote the branches of $(C,z)$ and \mbox{$J =
  J(C,T^\ast)$}.   Then
\[
f \in J \Leftrightarrow (f, Q_j) \ge 2\delta(Q_j) + \sum_{i\not= j} (Q_i,
Q_j) + \sum_{q \in T^\ast \cap Q_j} \mt \,(Q_{j,(q)}, q),
\]
where \mbox{$T^\ast \!\cap Q_j = \{q \in T^\ast \!\mid q \in Q_{j,(q)}\}$} and
$(f,g)$ denotes the intersection multiplicity.
\end{lemma}

\begin{proof}
As shown before, the multiplicity of the strict transform at \mbox{$q\in
T^\ast$} of a generic element
\mbox{$g\in J$} fulfills
\mbox{$\mt \,(g_{(q)},q)=\mt \,(X,q)=m_q$}. In particular, we obtain for each
branch $Q_j$ of $(C,z)$ 
\begin{eqnarray*}
\left(f,Q_j \right) \;  \geq \; \left( g,Q_j \right) & \ge &
\sum_{q\in T^\ast \cap Q_j} m_q \cdot \mt \,(Q_{j,(q)},q)\\
& = &  \sum_{q\in T^\ast \cap Q_j} \sum_{i=1}^r  \mt \,(Q_{i,(q)},q) \cdot
\mt \,(Q_{j,(q)},q)\\ 
& = & 2\delta(Q_j) + \sum_{i\neq j} (Q_i,Q_j) + \sum_{q\in T^\ast \cap Q_j}
\mt \,(Q_{j,(q)},q)\;=:\; \alpha_j\,. 
\end{eqnarray*}
Hence, we have the inclusion 
\mbox{$ J \subset J_1 := \bigcap_{j=1}^r \left\{ f\in \ko_{S,z} \mid
    (f,Q_j)\geq \alpha_j \right\}.$}
We can consider both as ideals in $\ko_{C,z}$ and have to show that
\mbox{$\dim_K (\ko_{C,z}/J_1) = \deg (X)$}. To do so, let
$n$ denote the injection 
$$ n: \ko_{C,z} \hookrightarrow \prod_{j=1}^r \C\{t_j\} =: \overline {\ko} $$
induced by a parametrization of $(C,z)$ and consider the image
\mbox{$n(J_1)\subset \overline{\ko}$}. For an element \mbox{$f \in
\overline{\ko}$} the conditions on the intersection multiplicities
\mbox{$(f,Q_j)$} read as \mbox{$f \in
\prod_{j=1}^r t_j^{\alpha_j}\cdot \C\{t_j\}$}. Hence
\begin{eqnarray*}
\dim_K\left(\ko_{C,z}\big/J_1\right) & \geq & \dim_K\big(\overline{\ko}\big/
\textstyle{\prod}_{j=1}^r \,t_j^{\alpha_j}\cdot \C\{t_j\}\big) -
\dim_K \big(
\overline{\ko}\big/ \ko_{C,z}\big)\\
& = & \sum_{j=1}^r \alpha_j - \delta(C,z)\; =\; \deg(X) .\qquad \qquad \qquad
\qquad \;\;\qed
\end{eqnarray*}
\renewcommand{\qed}{}\end{proof}

\begin{definition}\label{1.9}

Two generalized singularities \mbox{$X_0, X_1 \in \kg \ks$}, centred at $z$,
are called {\bf isomorphic}, $X_0 \cong X_1$, if they are isomorphic as
subschemes of $S$.   $X_0$ and $X_1$ are called {\bf equivalent}, \mbox{$X_0
\sim X_1$}, if there exist germs (respectively multigerms, if the $X_i$ are
reducible) \mbox{$(C_0, z)$} defining $X_0$ and \mbox{$(C_1, z)$}
defining $X_1$, 
and a $T^\ast \!$--equimultiple family over some (reduced) open connected
subset $T$ of $\A^1$ having $(C_0,T^\ast_0)$ and $(C_1,T^\ast_1)$ as
fibres. Here, by a {\bf $T^\ast \!$--equimultiple family}  over a (reduced)
algebraic $k$--scheme $T$ we denote a flat family with section $\sigma$,
\arraycolsep0.1cm
\renewcommand{\arraystretch}{1.1}
$$
\begin{array}{rrl}
\makebox[0.7cm]{$\kc$} & \hookrightarrow & S\times T\\
\makebox[0.3cm]{$\sigma \nwarrow$} & & \swarrow \\
& \makebox[0.1cm]{$T$ ,}
\end{array}
$$
of reduced plane curve singularities \mbox{$(\kc_t,\sigma(t))\subset S=S\times
\{t \}$} which admits a simultaneous, embedded
resolution, together with sections $\sigma_q$ through infinitely near points,
defining a family $\kt^\ast$ of trees $T^\ast_t$, \mbox{$T^\ast(\kc_t,
  \sigma(t)) \subset T^\ast_t \subset T(\kc_t,\sigma(t))$}
such that the total transform of $\kc$ is equimultiple along $\sigma_q$,
\mbox{$\sigma_q(t)\in T^\ast_t$}. 
We denote such a family by {\bf $(\kc,\kt^\ast)$}. 
\end{definition}

Since any \mbox{$X \in \kg\ks$} is defined by a generic element in \mbox{$J_X
  \subset \ko_{S,z}$}, isomorphic schemes $X_0, X_1$ are defined by isomorphic
  germs which can be connected by a family of isomorphisms.   Hence, \mbox{$X_0
  \cong X_1$} implies \mbox{$X_0 \sim X_1$}.

\begin{definition}\label{1.10}
  Let \mbox{$X \in \kg\ks$} with centre $z$ and \mbox{$L = (L,z)$} be a smooth
  (mul\-ti)germ. 
  Define \mbox{$X \cap L$} to be the scheme--theoretic intersection.   Set
  \begin{align*}
    T^\ast_X \cap L & := \{q \in T^\ast_X \mid q \in L_{(q)}\},\\ 
    J_X : L & :=\{f \in \ko_{S,z} \mid f L \in J_X\},\\ 
    X : L & := Z(J_X :L).
  \end{align*}

We call \mbox{$X : L$} the {\bf reduction of $\bf X$} by $L$.
\end{definition}

\begin{proposition}\label{1.11}
  Let \mbox{$X=X(C,T^\ast) \in \kg\ks$} and \mbox{$L \subset
  S$} be smooth at $z$, the centre of $X$.

  \begin{enumerate}
  \item[(i)] The reduction \mbox{$X :L$} is a generalized singularity centred
    at $z$ and its tree \mbox{$T^\ast_{X:L} = T^\ast\!:L$} is a subtree of
    $T^\ast$.

\item[(ii)] \mbox{$\mt \,(X,q)-1 \le \mt \,(X:L,q) \le \mt \,(X,q)$} for all
  \mbox{$q \in T^\ast$}.

\item[(iii)] 
\mbox{$\deg(X:L)  =   \deg(X) - \deg(X \cap L)$}, \mbox{$\;\deg(X \cap L)  =
  \! \underset{q \in T^\ast \cap L}{\sum} \mt \,(X,q)$}. 

\item[(iv)] Let $L$ be a line in $\P^2$, then there exists an exact sequence of
ideal sheaves on $\P^2$ 
\[
0 \lra  \kj_{X:L/\P^2}(d\!-\!1) \stackrel{\cdot L}{\lra} \kj_{X/\P^2}(d) \lra 
\kj_{X \cap L/L}(d)
\lra 0.
\]
  \end{enumerate}
\end{proposition}

\begin{proof}
(iv) is obvious and (iii) follows from (iv), respectively the fact that 
$$\deg (X\cap L) = \mt \,(C\cap L,z) = \!\sum_{q\in T^\ast \cap L} \!\mt
\,(X,q).$$ 
(i) will be proved by induction on $\deg(X)$.  Again, we may begin with
\mbox{$\deg(X) = 0$}, which implies \mbox{$X \!= \emptyset$}, \mbox{$T^\ast \!=
  \emptyset$} and
\mbox{$X : L = \emptyset$}, \mbox{$T^\ast \!: L = \emptyset$}.

If \mbox{$\deg(X) > 0$} then \mbox{$z \in T^\ast$} and we consider the
blowing--up \mbox{$\pi : \widehat{S} \lra S$} of $z$. By Lemma \ref{1.6} the
strict transform $C^\ast$ of $C$ fulfills \mbox{$\deg(X(C^\ast \!,
  T^\ast\backslash\{z\})) < \deg(X)$},
hence by the induction assumption there exists a (multi)germ $D^{\ast}$ and
a subtree \mbox{$T_{D^\ast} \subset T^\ast \backslash\{z\}$} such that
\begin{gather*}
J(C^{\ast}, T^\ast \backslash\{z\}) : L^\ast = J(D^{\ast}, T_{D^\ast})\\
T^\ast(D^{\ast}) \subset T_{D^\ast} \subset T(D^{\ast}),
\end{gather*}
moreover, we can choose $D^{\ast}$ generically in \mbox{$J(C^\ast, T^\ast
  \backslash\{z\}) : L^\ast$} (such that
$\mt \,(D^\ast)$ and the intersection multiplicity $(D^\ast \!, E)$ are
minimal).

Blowing down $D^{\ast}$ we obtain a germ $D$ at $z$. Let \mbox{$m:=\mt
  \,(C,z)$}.

\underline{Case 1}:  $\mt \,(D,z) = m-a < m$.

Define the germ
$C^\prime$ at $z$ by
\mbox{$C^\prime := D \cdot L_1 \cdot \ldots \cdot L_{a-1}$},
where \mbox{$L_1, \dots, L_{a-1}$} are smooth germs with generic tangent
directions at $z$.   Then
\[
\begin{array}{lcl}
f \in J(C,T^\ast) : L & \!\Leftrightarrow & \!\mt \,(f,z) \ge m\!-\!1,\;L^\ast
\hat{f} \!:\! (m\!-\!1) E \in J(C^{\ast}, T^\ast \backslash\{z\})\\
& \!\Leftrightarrow & \!\mt \,(f,z) \ge m\!-\!1,\; \hat{f}\! :\! (m\!-\!1) E
\in 
J(D^{\ast},T_{D^\ast})\\
& \!\Leftrightarrow & \!\mt \,(f,z) \ge m\!-\!1, \;\hat{f}\! :\! (m\!-\!1) E
\in 
J(D^{\ast} \!\!\cdot
L^\ast_1 \!\cdot \ldots \cdot L^\ast_{a-1}, T_{D^\ast})\\
& \!\Leftrightarrow & \!f \in J (C^\prime, T^{\prime\ast})
\end{array}
\]
with \mbox{$T^{\prime\ast} := T_{D^\ast} \cup \{z\}$}. 

\underline{Case 2}: \mbox{$\mt \,(D,z) = m$}.

By the induction assumption, there exists a (multi)germ
$\bar{D}^\ast$ and a subtree \mbox{$T_{\bar{D}^\ast} \subset T_{D^\ast} \subset
    T^\ast \backslash\{z\}$} such that
\begin{gather*}
J(D^{\ast}, T_{D^\ast}) : E = J(\bar{D}^\ast, T_{\bar{D}^\ast})\\
T^\ast(\bar{D}^\ast) \subset T_{\bar{D}^\ast} \subset T(\bar{D}^\ast).  
\end{gather*}

We choose $\bar{D}^\ast$ generically in the ideal \mbox{$J(D^\ast, T_{D^\ast})
  : E$}. Since \mbox{$D^{\ast}\! \in
  J(\bar{D}^\ast\!, T_{\bar{D}^\ast})$},
\mbox{$p:= (D^{\ast}\!, E)\! -\! (\bar{D}^{\ast}\!, E) \ge 0$}, and we 
define \mbox{$C^\prime :=\bar{D} \cdot L_1 \cdot \ldots \cdot L_p$},
where \mbox{$L_1, \dots, L_p$} denote generic smooth germs at $z$ and
$\bar{D}$ 
is the blowing--down of $\bar{D}^\ast$.   Again
\[
f \in J(C,T^\ast) :L \;\Leftrightarrow \; \mt \,(f,z) \ge m\!-\!1, \;L^\ast
\hat{f} : (m\!-\!1) E \in J (C^{\ast}, T^\ast \backslash \{z\}).
\]

Assuming \mbox{$\mt \,(f,z) = m\!-\!1$} and \mbox{$L^\ast \hat{f} :(m\!-\!1) E
  \in
  J(C^{\ast}, T^\ast \backslash\{z\})$} we would have \mbox{$f^\ast \!\in
  J(C^{\ast}\!,  T^\ast \backslash\{z\})
  : L^\ast = J(D^{\ast}\!,T_{D^\ast})$}, hence \mbox{$m\!-\!1 = (f^\ast \!, E)
  \ge (D^{\ast} \!, E) = m$}. Thus
\[
\begin{array}[l]{lcl}
f \in J(C,T^\ast) : L & \Leftrightarrow & \mt \,(f,z) \ge m,
\;\: \hat{f} : mE\in J(D^{\ast}, T_{D^\ast}) :\!E = J(\bar{D}^\ast,
T_{\bar{D}^\ast})\\
& \Leftrightarrow & \mt \,(f,z) \ge m,\;\: \hat{f} : mE \in J(\bar{D}^\ast
\cdot
L^\ast_1 \cdot \ldots \cdot L^\ast_p, T_{\bar{D}^\ast})\\
& \Leftrightarrow & f \in J(C^\prime,T^{\prime\ast}\!:= T_{\bar{D}^\ast}
\cup \{z\}).  
\end{array}
\]
Note that in this case we have \mbox{$\mt \,(C^\prime,z)=\mt \,(C,z)$}, while
in the
first case we had \mbox{$\mt \,(C^\prime,z)=\mt \,(C,z)-1$}. This implies (ii).
\end{proof}

\begin{examples} 
%
%
\begin{enumerate}
\itemsep0.1cm
\item Let \mbox{$(C,z)$} be a node, \mbox{$C=y^2-x^2$} with
respect to local coordinates \mbox{$(x,y)$} at $z$,
\mbox{$T^\ast=T^\ast(C,z)=\{z\}$}, then for each $L$ the reduction \mbox{$X:L$}
is the
generalized singularity given by a smooth
germ at $z$ and the tree \mbox{$T^\ast_{X:L}=T^\ast=\{z\}$}.
\item  In the case of an $A_{2k-1}$--singularity (\mbox{$k\geq 2$})
  \mbox{$C=y^2-x^{2k}$} with the tree of essential points
\mbox{$T^\ast=\{z=q_0,\dots,q_{k-1}\}$} we have
  \mbox{$J_X=\langle
  y^2,yx^k,x^{2k}\rangle $}. If \mbox{$L=y$}, then \mbox{$J_{X:L}=\langle
  y,x^k\rangle$} and \mbox{$X:L$} is the generalized singularity given by the
smooth germ
\mbox{$y-x^k$} and the tree \mbox{$T^\ast_{X:L}=\{z=q_0,\dots ,q_{k-1}\}$}. 
On the other hand, if $L=x$, then  \mbox{$J_{X:L}=\langle
  y^2,yx^{k-1},x^{2k-1}\rangle$} and \mbox{$X:L$} is the generalized
singularity scheme given by \mbox{$(y-x^k)(y+x^{k-1})$}
and the tree
\mbox{$T^\ast_{X:L}=T^\ast $}.
\item For an $A_{2k}$--singularity \mbox{$C=y^2-x^{2k+1}$} with the tree of
  essential points
\mbox{$T^\ast=\{z=q_0,\dots,q_{k+1}\}$} we have \mbox{$J_X=\langle
  y^2,yx^{k+1},x^{2k+1}\rangle $}. If \mbox{$L=y$}, then \mbox{$J_{X:L}=\langle
  y,x^{k+1} \rangle$} and \mbox{$X:L$} is the generalized singularity given by
the smooth germ
\mbox{$y-x^{k+1}$} and the tree
\mbox{$T^\ast_{X:L}=\{z=q_0,\dots,q_k\}$}. On the other
hand, if $L=x$, then  \mbox{$J_{X:L}=\langle y^2,yx^k,x^{2k}\rangle $} and
\mbox{$X:L$} is
the singularity scheme given by an $A_{2k-1}$--singularity. 
\end{enumerate}
\end{examples}

\begin{definition}\label{1.12}
  Denote by \mbox{$\kg\ks_1 \subset \kg\ks$} the subclass of such $X$ defined
  by germs $(C,z)$ with all branches smooth.
\end{definition}

\smallskip
\begin{lemma}\label{1.13}
  The class $\kg\ks_1$ is closed with respect to the equivalence relation
  $\sim$ and with respect to reduction by  $L$.
\end{lemma}
\begin{proof}
The first statement is obvious, the second is a consequence of
the proof of \ref{1.11}.
\end{proof}

\begin{lemma}\label{1.14}
  Let \mbox{$X = X(C,T^\ast) \in \kg\ks$} be non--empty and $L$  smooth at
  $z$.

  \begin{enumerate}
  \item[(i)] There exists a branch $Q$ of $(C,z)$ such that \mbox{$T^\ast \cap
      L \subset T^\ast \cap Q$}.

\item[(ii)] If $Q$ is a non--singular branch of $(C,z)$ and if \mbox{$M \subset
  T^\ast$}
  is a connected subtree with \mbox{$T^\ast\!\cap L \subset M \subset
  T^\ast\!\cap Q$} then there exists \mbox{$X_1 = X(C_1, T^\ast_1) \in
  \kg\ks$}, 
  \mbox{$X_1 \cong X$},
  $Q_1$ a smooth branch of \mbox{$(C_1, z) \cong (C,z)$}, \mbox{$M \cong M_1
  \subset Q_1 \cap T^\ast_1$} such that \mbox{$T^\ast_1 \cap L = M_1$}.  
  \end{enumerate}
\end{lemma}

\begin{proof}
(i) is obvious. For the proof of (ii), we may assume that \mbox{$z\in L\cap
T^\ast$}. We choose coordinates $(x,y)$ at z such that \mbox{$L=y$} and the
non--singular branch $Q$ is given by \mbox{$y-\sum_{i=1}^\infty \alpha_i
x^i$}. Let \mbox{$N:=\#(M)$} and $f$ be the power series defining
\mbox{$(C,z)$}. The germs 
$$ C_t\,=\,f\left(x\,,\,y+t\cdot \!\left(\textstyle{\sum\limits_{i=1}^{N-1}}
\alpha_i x^i+\beta x^N\right)\right)\,, \;\;\beta \in K \mbox{ generic,}$$
define an equianalytic family such that \mbox{$(C_0,z)=(C,z)$} and
\mbox{$(C_1,z)$} has a branch $Q_1$ given by \mbox{$y-\sum_{i\geq N}
\tilde{\alpha}_i x^i$}. Especially for the corresponding trees
\mbox{$T^\ast_1\cong T^\ast$}, respectively \mbox{$M_1\cong M$}, we have
\mbox{$M_1 \subset Q_1 \cap T^\ast_1$} and, since $\beta$ was chosen
generically,  
$$ \#(T^\ast_1 \cap L )= \#(T^\ast_1\cap L \cap Q_1) = N = \#(M_1), $$
hence \mbox{$T^\ast_1\cap L = M_1$}.
\end{proof}

\pagebreak[3]
\begin{lemma}\label{1.15}
  Let \mbox{$X = X(C,T^\ast) \in \kg\ks$}, $Q$ be a smooth branch of $(C,z)$
  and $L$ be smooth at $z$.

  \begin{enumerate}
  \item[(i)] If \mbox{$T^\ast \cap Q \subset T^\ast \cap L$} then
    \mbox{$\mt \, (X:L) = \mt \, X -1$}.
\item[(ii)] If \mbox{$T^\ast \cap Q = T^\ast \cap L$} then \mbox{$X :L$} is
  defined by the germ $(C^\prime \!,z)$ and the tree \mbox{$T^\ast\!\cap
T(C^\prime \!,z)$}, where $C^\prime$ is a factor of $C$
such that \mbox{$C = C^\prime Q$}. 
  \end{enumerate}
\end{lemma}

\begin{proof}
Let \mbox{$C=C^\prime\cdot Q$} at $z$. Then, obviously
\mbox{$C^\prime \in J_X:L$}, hence the
multiplicity of a generic element is at most \mbox{$\mt \,X-1$}. Thus, (i)
follows from \ref{1.11} (ii).\\
 If \mbox{$T^\ast \!\cap Q = T^\ast \!\cap L$},
then we have for each
\mbox{$q\in T^\ast$} 
$$ \mt \,(\widehat{C}^\prime_{(q)},q)=
\mt \,(\widehat{C}_{(q)},q)-\mt \,(\widehat{L}_{(q)},q),$$  
which implies (ii).
\end{proof}

\begin{lemma}\label{1.16}
  Let \mbox{$X = X(C,T^\ast) \in \kg\ks$} and $(L,z)$ be a smooth germ.  Then
  we have   \mbox{$\deg((X:L) \cap L) \le \deg(X \cap L)$}. Moreover,

  \begin{enumerate}
  \item[(i)] if \mbox{$\deg((X:L) \cap L) = \deg (X \cap L)$} then
    \mbox{$\mt \,(X:L) = \mt \,X$};
\item[(ii)] if \mbox{$\deg((X:L) \cap L)) < \deg(X \cap L)$} then
either \mbox{$\mt (X:L) < \mt \, X $} or 
the defining germ $(C',z)$ of \mbox{$X:L$} has a branch $Q'$
  satisfying 
\mbox{$ T^\ast_{X:L} \cap Q' \subset   T^\ast_{X:L} \cap L$}. 
In any case, \mbox{$\mt (X:L^2) \le \mt \, X -1$}.
  \end{enumerate}
\end{lemma}

\begin{proof}
  By Prop.\
\ref{1.11} (ii), \mbox{$\mt \,(X:L,q) \le \mt \,(X,q)$} for all \mbox{$q
  \in T^\ast$}, hence the inequality.   Therefore, \mbox{$\deg((X:L) \cap L) =
  \deg(X \cap L)$} implies \mbox{$\mt \,(X:L,q) = \mt \,(X,q)$} for 
any \mbox{$q \in T^\ast\!
  \cap L$}, in particular for \mbox{$q = z$}.   This implies (i).\\
Let \mbox{$\deg((X:L) \cap L) < \deg(X \cap L)$} and \mbox{$T^\ast\!
  \cap L=\{z\!=\!q_0,q_1,\dots,q_\ell\}$}.
Recall the construction of \mbox{$X:L=X(C^\prime,T^{\prime\ast})$} in
  the proof of Proposition \ref{1.11} (i). In \mbox{Case 1} we had 
$$ \mt \,(X:L)=\mt \,(C^\prime \!,z)=\mt \,(C,z)-1=\mt \,X-1. $$
In Case 2, $C^\prime$ was given as 
$$ C^\prime=\bar{D}\cdot L_1 \cdot \ldots \cdot L_p \:\mbox{ with
  }\:p=(C^\ast \!,E)-(\bar{D}^\ast \!,E).$$ 
 Assume that there is no branch $Q'$ of
  $C^\prime$ such that \mbox{$T^{\prime\ast} \!\cap Q' \subset
  T^{\prime\ast} \!\cap L$}, in particular,
  \mbox{$p=0$}. Then \mbox{$\bar{D}^\ast\!=C^{\prime\ast}$}, the
  strict transform of $C^\prime$, and
  $$\mt \,(C_{(q_0)})=(C^\ast \!,E)=(C^{\prime\ast} \!,E)=
  \mt \,(C^\prime_{(q_0)}).$$
 On the other hand, the intersection multiplicity
  \mbox{$(C^{\prime\ast} \!,E)$} is just the sum of all
  \mbox{$\mt \,(C^\prime_{(q)})$}
  with \mbox{$q\in T(C^\prime)\cap E_{(q)}$}. By the above assumption all those
  $q$ are essential for $C^\prime$, which implies 
$$ \sum_{q\in T^\ast(C^\prime)\cap E_{(q)}} \mt \,(C^\prime_{(q)})=
(C^{\prime\ast} \!,E)
= (C^\ast \!,E)= \sum_{q\in T(C)\cap E_{(q)}} \mt \,(C_{(q)}).$$
Since  \mbox{$T^\ast(C^\prime)\subset T^\ast(C)$}, it follows
\mbox{$\mt \,(C^\prime_{(q)})=\mt \,(C_{(q)})$} for all \mbox{$q\in
  T^\ast(C^\prime)\cap
    E_{(q)}$}, especially \mbox{$\mt \,(C^\prime_{(q_1)})=\mt \,(C_{(q_1)})$}.
By induction, we obtain \mbox{$\mt \,(C^\prime_{(q_i)})=\mt \,(C_{(q_i)})$} for
each \mbox{$i\in \{0,\dots\!\,,\ell\}$}, which is 
impossible (cf.~Proposition \ref{1.11} (iii)). 
\end{proof}

\begin{definition}
Given a $T^\ast\!$--equimultiple family of plane curve singularities
\mbox{$(\kc,\kt^\ast)$} over a reduced algebraic $K$--scheme $T$ (as defined in
\ref{1.9}), we define 
$$
\mathfrak{J}(\kc,\kt^\ast):=\left\{ f\!\in\! \ko_{S\times T,\sigma(T)} \,\big|
  \, 
\mt \,(\hat{f}_{\sigma_q(t)},\sigma_q(t))\geq
\mt \,(\hat{\kc}_{\sigma_q(t)},\sigma_q(t) ) \mbox{ for } q \!\in \!T^\ast_t,\,
t\!\in \!T\right\}
$$
and
$$\kx(\kc,\kt^\ast):= \ko_{S\times T,\sigma(T)}\big/\mathfrak{J}
(\kc,\kt^\ast).$$
A flat family $\kx$ of fat points in \mbox{$S\times T$} is called a {\bf family
of generalized singularity schemes}, if \mbox{$\kx =\kx(\kc,\kt^\ast)$} for
some $T^\ast\!$--equimultiple family \mbox{$(\kc,\kt^\ast)$}.
\end{definition}

Since we consider only
reduced base spaces $T$, then flatness just means that the total length is
constant, which holds for a family  \mbox{$\kx(\kc,\kt^\ast)$} by Lemma
\ref{1.6}. It is easily seen that the functor 
$$\underline{\kg\ks} \,: \:T \mapsto \{\mbox{families of generalized
singularity
schemes over }T \}$$
is representable by a locally closed subscheme $GS$ of the punctual
Hilbert scheme of $S$.

\begin{proposition}\label{1.17a}
  Let \mbox{$X \in \kg\ks$}, $L$ be smooth and \mbox{$Y = X:L$}.   For almost
  all \mbox{$Y^\prime \sim Y$} satisfying
  \mbox{$\deg(Y^\prime \cap L) = \deg(Y \cap L)$} there exists
  a generalized singularity scheme \mbox{$X^\prime \sim X$} such that
  \mbox{$\deg(X^\prime \cap L) = \deg(X \cap 
  L)$} and \mbox{$Y^\prime = X^\prime :L$}.
\end{proposition}

\begin{proof}
\footnote{We should like to thank I.~Tyomkin for an idea leading to the present
  proof.}
Let \mbox{$\kx =\kx(\kc,\kt^\ast)$} be a family of generalized singularity
schemes over the reduced base space $T$, \mbox{$t\in T$}. The construction of
\mbox{$X:L$} given in the proof of \ref{1.11} shows that we can simultaneously
reduce the fibres of $\kx$ by $L$. Hence we have a natural transformation
$$ \rho_L: \underline{\kg\ks}\lra \underline{\kg\ks}\:,\; \kx \mapsto \kx :L $$
inducing a morphism \mbox{$\rho_L : GS \lra GS$}. Notice that two generalized
singularity schemes $X_1$,$X_2$ are equivalent if and only if they are in the
same connected component of $GS$. Therefore, to prove the proposition, it is
enough to show that the restriction 
\mbox{$ \rho_{L,X}: GS_{L,X} \lra GS_{L,Y}$}
of $\rho_L$ to the connected component $GS_{L,X}$ of
\mbox{$\{X^\prime \in GS \mid \deg(X^\prime \!\cap L)= \deg(X\cap L)\}$}
containing $X$ is
dominant. But this follows immediately from the fact that the dimension of the
fibre \mbox{$ \rho_{L,X}^{-1}(Y)$} is just
\mbox{$
\#(T^\ast_X)-\left(\#(T^\ast_Y)+\#(T^\ast_X \cap L) -\#(T^\ast_Y \cap L)\right)
= \dim (GS_{L,X}) - \dim (GS_{L,Y}).$} 
\end{proof}

In the following, we shall introduce the second basic operation on generalized
singularities, the extension.   For this, it is convenient to work with the
field $K\{\{x\}\}$ of fractional power series
\[
\sum^\infty_{i=0} \alpha_i x^{i/n}, \alpha_i \in K, n \in \N.
\]

Any germ $(C,z)$ of a reduced curve singularity may be given, with respect to
suitable local coordinates $x,y$, as
\[
C = \prod^m_{i=1} (y - \xi_i(x)),\;\; m = \mt \,(C,z),\;\; \xi_1, \dots, \xi_m
\in K\{\{x\}\}.
\]
Moreover, if $(C,z)$ is irreducible and $t=x^{1/m}$, then 
$$C = \prod^m_{i=1} (y - \xi(\eta^it)),\;\; \xi \in K[[t]],$$
with $\eta$ a primitive $m$--th root of unity.
We define the intersection multiplicity of two fractional power series $\xi_i,
\xi_j\in K\{\{x\}\} $ to be
\[
(\xi_i, \xi_j) := \max \,\{\rho \in \Q \mid x^\rho \mbox{ divides } \xi_i(x) -
\xi_j(x)\}.
\]

\begin{lemma}\label{1.17}
  Let \mbox{$X = X(C,T^\ast(C,z)) \in \ks$} be a singularity scheme  and
  $(C,z)$ given as above.   Then
\[
\deg(X) = \sum_{1\le i < j \le m} (\xi_i, \xi_j) + \sum^m_{i=1} \underset{j
  \not=i}{\max} \,(\xi_i, \xi_j) + \frac{m-r}{2},
\]
where \mbox{$m = \mt \,(C,z)$} and $r$ is the number of branches of $(C,z)$.
\end{lemma}

\begin{proof}
 It is well--known that the intersection multiplicity at $z$ of
the polar 
$P_q(C)$ (\mbox{$q=(0\!:\!1\!:\!0)$}) given by the power series
\[
\frac{\partial C}{\partial y} = \sum_{i=1}^m \prod_{j\neq i} (y-\xi_j(x)) 
\]
and the curve C fulfills
\mbox{$ 
\sum_{i\neq j} (\xi_i,\xi_j)=\mt \,(P_q(C)\cap C,z)=2\delta(C,z)+m-r$}.\\ 
Hence, by Lemma \ref{1.6}, it suffices to show
\[
\sum_{q\in T^\ast(C)} m_q=\sum_{i=1}^m \underset{j\neq i}{\max}
\,(\xi_i,\xi_j)+m-r. 
\]
In the case of an irreducible germ  $(C,z)$ it follows from the above
description of the $\xi_i$ that the numbers $(\xi_i,\xi_j)$ do only depend on
the characteristic terms of the Puiseux expansion, and the statement is an
immediate consequence of the algorithm to compute the multiplicity sequence
from the Puiseux pairs (cf.~\cite{BK}).\\
In the case of a reducible germ $(C,z)$, we have to investigate, additionally,
the case of two branches \mbox{$Q_k=\prod_{i=1}^{m_k} (y-\xi^{(k)}(\eta_k^i
t))$} 
\mbox{($k\!\in\!\{1,2\}$)} such that \mbox{$T(Q_1)\cap T(Q_2)$} contains a
non--essential point $q$ of $T(Q_1)$ and for all branches \mbox{$Q\neq Q_1$}
of \mbox{$(C,z)$} and all successors $\hat{q}$ of $q$ in \mbox{$T(Q_1)$} we
have \mbox{$\hat{q}\not\in T(Q)$}.  In this case, obviously,
\mbox{$m_2=M
m_1$}, \mbox{$M\in \N$}, and we can assume the maximum intersection
multiplicity of the fractional power series \mbox{$\xi^{(1)}_i
(x)=\xi^{(1)}(\eta^{Mi}x^{1/m_1})$} ($\eta $
a primitive $m_2$--th root of unity, \mbox{$i\in \{1,\dots,m_1\}$}) with any
other fractional power series in
the equation of $(C,z)$ to be realized by \mbox{$\xi^{(2)}_{Mi}
(x)=\xi^{(2)}(\eta^{Mi}x^{1/m_2})$}. Then we have
\begin{eqnarray*}
\sum_{q\in T^\ast(C)} m_q(Q_1) & = & \frac{1}{M} \sum_{q\in T^\ast(C)}
m_q(Q_1)m_q(Q_2) - \sum_{q\in T^\ast(C)} m_q(Q_1)(m_q(Q_1)-1) \\
&=& \frac{1}{M} \cdot \mt \,(Q_1\cap Q_2,z)-2\cdot \delta (Q_1)\\
& = & \sum_{i=1}^{m_1} \sum_{j=1}^{m_1}
(\xi_i^{(1)},\xi_{Mj}^{(2)}) - \sum_{i\neq j} (\xi_i^{(1)},\xi_j^{(1)}) + m_1
-1,
\end{eqnarray*}
and the statement follows from the fact that for \mbox{$i\neq j$} the
intersection multiplicities  \mbox{$(\xi_i^{(1)},\xi_j^{(1)})$} and \mbox{$
(\xi_i^{(1)},\xi_{Mj}^{(2)})$} coincide. 
\end{proof}

\begin{lemma}\label{1.18}
  Let \mbox{$X = X(C,T^\ast) \in \kg\ks$}, $L$  be smooth at $z$ and \mbox{$q
  \in T^\ast \!\cap L\backslash \{z\}$}.   Let
\[
C = \prod^n_{i=1} (y - \xi_i(x)) \prod^m_{i=n+1} (y-\xi_i(x)),\;\; m =
\mt \,(C,z),\;\; \xi_i \in K\{\{x\}\},
\]
be decomposed so that \mbox{$\xi_1, \dots, \xi_n$} are all fractional power
series belonging to branches $Q$ of $(C,z)$ with $q \in Q_{(q)}$.   Then there
exists an integer \mbox{$k \ge 0$} such that
\begin{align*}
  k < (\xi_i, \xi_j) & \mbox{ for } 1 \le i < j \le n,\\
k \ge(\xi_i, \xi_j) & \mbox{ for } 1 \le i \le n < j \le m.
\end{align*}

More precisely, if \mbox{$L=y$}, then
\[
\xi_i(x) = \sum_{\rho \ge 0} \alpha^{(i)}_\rho x^\rho,\quad \alpha^{(i)}_\rho
\in K,\; \rho \in \Q
\]
belongs to $Q$ with \mbox{$q \in Q_{(q)}$} if and only if
\mbox{$\alpha^{(i)}_\rho \!= 0$} for \mbox{$\rho \le k$}. 
\end{lemma}

\begin{proof}
Let $L=y$ and \mbox{$T^\ast \!\cap
L = \{z\!=\!q_0,q_1,\dots,q_\ell\}$}. Moreover, let the branch $Q$ be given by
\[
Q=y^p+a_1(x)y^{p-1}+\dots +a_p(x)= \prod_{i=1}^p\left(y-\xi (\eta^i
x^{1/p})\right),
\]
where \mbox{$\xi(t)=\sum_{j=0}^\infty \alpha_j t^{j} \in
K[[t]]$}. To prove the lemma, it is
sufficient to show that for each \mbox{$k\in \{1,\dots \!\,,\ell \}$}
\[
q_k\in Q_{(q_k)} \;\;\Leftrightarrow \;\; \alpha_j=0 \,\mbox{ for each }\,
j\leq k\!\cdot \! p.
\]
We proceed by induction on the length $k+1$ of the tree \mbox{$\{q_0,\dots,
q_k\}$}.
Obviously, \mbox{$z=q_0\in Q$} if and only if \mbox{$a_p(0)=0$},
that is, if and only if \mbox{$\alpha_0=0$}. Furthermore, the total transform
of
$Q$ at $q_1$ reads as \mbox{$\widehat{Q}_{(q_1)}=\prod_{i=1}^p (uv-\xi(\eta^i
u^{1/p}))$}, hence \mbox{$q_1\in Q_{(q_1)}$} if and only if
\mbox{$\xi(t)=\sum_{j=p+1}^\infty \alpha_j t^{j} $}. Then $Q_{(q_1)}$ has the
equation 
\[ 
Q_{(q_1)} = \prod_{i=1}^p \left(v-\tilde{\xi}(\eta^iu^{1/p})\right)
\]
at $q_1$, where \mbox{$\tilde{\xi}(t)=\sum_{j=1}^\infty \alpha_{j+p}
t^{j}$}, and we complete the proof  by applying the induction
hypothesis to $Q_{(q_1)}$ and the tree \mbox{$\{q_1,\dots , q_k\}\subset
\left(T^\ast\backslash\{z\}\right)\cap L$}.
\end{proof}

\begin{definition}\label{1.19}
  Using the notations and hypotheses of Lemma \ref{1.18}, let
\[
\xi_i(x) = \sum_{\rho > k} \alpha_\rho^{(i)} x^\rho,\; i = 1, \dots, n.
\]

Define a germ \mbox{$(C(q), z)$} by
\[
C(q) := \prod^n_{i=1} (y - x \xi_i(x)) \prod^m_{i = n+1} (y - \xi_i(x)).
\]

Call $C(q)$ the {\bf extension of} $\bf C$ at $q$.
\end{definition}

\begin{lemma}\label{1.20}
  The tree $T^\ast(C(q))$ of essential points of $C(q)$ has the following
  structure:  insert in $T^\ast(C)$ a new point $q^\prime$ between $q$ and its
  predecessor $\bar{q}$.   Moreover, \mbox{$\mt \,(C(q)_{(p)},p) =
  \mt \,(C_{(p)},p)$}
  for all \mbox{$p \in  T^\ast(C(q))\backslash \{q^\prime\} = T^\ast(C)$} and
\[
\mt \,(C(q)_{(q^\prime)}, q^\prime) = \sum_{Q:\;q\in Q_{(q)}} \mt
\,(Q_{(\bar{q})}, \bar{q}).
\]
\end{lemma}

Any tree $T^\ast$ containing $T^\ast(C)$ becomes extended by this operation to
a tree $T^\ast(q)$.   We call $T^\ast(q)$ the {\bf extension of} $\bf T^\ast$
at $q$.

\begin{proof}
As in the proof of \ref{1.18}, let \mbox{$T^\ast\!\cap L=\{z\!=\!q_0,q_1,\dots
, q_\ell\}$}. An easy consideration shows that for \mbox{$q=q_k$} the strict
transform \mbox{$C(q)_{(q)}$} of the extension of $C$ at $q$ has the local
equation 
\[
\prod_{i=1}^n \big(v-\frac{1}{u^{k-1}}\xi_i(u)\big) =
\prod_{i=1}^n \big(v-\textstyle{\sum\limits_{\rho >0}} \alpha^{(i)}_{\rho +k}
u^{\rho +1}\big) 
\]
while \mbox{$C(q)_{(q_{k+1})}$} is given by
\[
\prod_{i=1}^n \big(v-\frac{1}{u^{k}}\xi_i(u)\big) =
\prod_{i=1}^n \big(v-\textstyle{\sum\limits_{\rho >0}} \alpha^{(i)}_{\rho +k}
u^{\rho }\big)
\]
which corresponds to the equation of $C_{(q)}$ at $q$. Hence, the structure of
\mbox{$T^\ast(C(q))$} can be described as in the lemma. Moreover, notice that
for \mbox{$i\leq k$}
$$ \mt \,(C(q)_{(q_i)},q_i)= n+ \!\sum_{Q:\;q\not\in Q_{(q)}}
\!\mt \,(Q_{(q_i)},q_i),$$
which implies the statement about the multiplicities.
\end{proof}

\begin{definition}\label{1.21}
  Let \mbox{$X = X(C,T^\ast) \in \kg\ks$}, $L$ be a smooth germ at $z$ and
  $q \in T^\ast \!\cap L \backslash \{z\}$.   We define the {\bf
  extension of} $\bf X$ at $q$ to be
\[
X(q) := X(C(q),\, T^\ast(q)) \in \kg\ks.
\]
\end{definition}

\begin{examples} 
\begin{enumerate}
\itemsep0.1cm
\item  Let \mbox{$(C,z)$} be an ordinary cusp, \mbox{$C=y^2-x^3$} with
respect to local coordinates \mbox{$(x,y)$} at $z$,
\mbox{$T^\ast=T^\ast(C,z)=\{z=q_0,q_1,q_2\}$} with (strict) multiplicities
\mbox{$m_z=2$}, \mbox{$m_{q_1}=m_{q_2}=1$}; we write \mbox{$T^\ast=
  \underset{2}{\ast}-\underset{1}{\ast}-\underset{1}{\ast}$}\,.\\
Moreover, let \mbox{$L=y$}, that is \mbox{$T^\ast \cap L=\{q_0,q_1\}$}. The
extension $X(q_1)$ is given by \mbox{$y^2-x^5$} and the tree
\mbox{$T^\ast(q_1)=
  \underset{2}{\ast}-\underset{2}{\ast}-\underset{1}{\ast}-
  \underset{1}{\ast}$}\,.

\item Let \mbox{$X\in \ks$} be given by \mbox{$C=(y^2-x^3)(y^2-x^5)$} and the
  tree of essential points
$$T^\ast=\underset{4}{\ast}- \underset{3}{\ast}
<\!\!\!
\renewcommand{\arraystretch}{0.7}
\begin{array}{l}
\overset{1}{\ast}-\overset{1}{\ast}\\
\underset{1}{\ast}
\end{array}.
$$
If \mbox{$L=y$}, then \mbox{$T^\ast\cap L=\{z=q_0,q_1,q_2\}$} and the extension
$X(q_1)$ is given by \mbox{$(y^2-x^5)(y^2-x^7)$} and the tree
 $$T^\ast(q_1)=\underset{4}{\ast}- \underset{4}{\ast}-\underset{3}{\ast}
<\!\!\!
\renewcommand{\arraystretch}{0.7}
\begin{array}{l}
\overset{1}{\ast}-\overset{1}{\ast}\\
\underset{1}{\ast}
\end{array},
$$
whence the extension
$X(q_2)$ is given by \mbox{$(y^2-x^3)(y^2-x^7)$} and  
$$T^\ast(q_2)=\underset{4}{\ast}-\underset{3}{\ast}
<\!\!\!
\renewcommand{\arraystretch}{0.7}
\begin{array}{l}
\overset{2}{\ast}-\overset{1}{\ast}-\overset{1}{\ast}\\
\underset{1}{\ast}
\end{array}.
$$
\end{enumerate}
\end{examples}

\begin{lemma}\label{1.22}
  With the notations of the preceding definition, assume that
  $$H^1(S,\kj_{X(q)/\P^2}(d)) = 0.$$ 
Then there exists an \mbox{$X^\prime \in \kg\ks$}, \mbox{$X^\prime \sim X$}
such that \mbox{$\deg(X^\prime \cap L)  = \deg(X \cap L)$} and
\mbox{$H^1(S,\kj_{X^\prime/\P^2}(d))  = 0$}.
\end{lemma}

\begin{proof}
Let \mbox{$(x,y)$} be coordinates in a neighbourhood of \mbox{$z=(0,0)$}, such
that \mbox{$L=y$} and $(C,z)$ is given as in Lemma \ref{1.18}. For
\mbox{$t\in \A^1$} define 
\[
C_t:=
\underbrace{\prod_{i=1}^n\big(y-((1\!-\!t)x+t)\cdot\xi_i(x)\big)}_{\textstyle
C^1_t} 
\:\cdot \underbrace{\prod_{i=n+1}^m\big(y-\xi_i(x)\big)}_{\textstyle C^2}.
\]  
For \mbox{$t\neq 0$} there is an obvious isomorphism \mbox{$\varphi_t:
(C^1_t,z)\stackrel{\cong}{\lra} (C^1_1,z)=(C^1\!,z)\subset (C,z)$} and for two
branches $Q^1_t$ (resp.~$Q^2$) of \mbox{$(C^1_t,z)$} (resp.~\mbox{$(C^2,z)$}
the intersection mul\-ti\-plicity at $z$ fulfills \mbox{$\mt \,(Q^1_t\cap
  Q^2,z)= \mt \,(Q^1_1\cap Q^2,z)$}, hence \mbox{$(C_t,z)\sim 
(C,z)$}. Moreover, for \mbox{$t\neq 0$} sufficiently small, $C_t$ has an
ordinary $n$--fold point at $z_t=(-t/(1\!-\!t),0)$.\\
Define $T^\ast_t$ as the union of $\{z_t\}$ with the tree $T^{2\ast}$
corresponding to $C^2$ and the tree induced by $\varphi_t$ from $T^{1\ast}$
(corresponding to $C^1$). $T^\ast_t$ is well--defined since
\mbox{$L\cap C_t\supset L\cap C$} for each \mbox{$t\in \A^1$}. Thus, we have
defined a family $\kx$ with fibres \mbox{$X_t=X(C_t,T^\ast_t)$} centred at the
multigerm \mbox{$\{z,z_t\}$}. Obviously $\kx$ is flat in \mbox{$t=0$} since for
small \mbox{$t\neq 0$} (by
\ref{1.6} and \ref{1.20}) 
$$\deg (X_t) = \deg (X)+\frac{n(n+1)}{2} = \deg (X(q))= \deg(X_0).$$
Hence, the family $\mathfrak{J}$ of ideals \mbox{$J_{X_t}=J(C_t,T^\ast_t)$} is
flat in \mbox{$t=0$}, which implies, by semicontinuity, the vanishing of
\mbox{$H^1(S,\kj_{X_t/\P^2}(d))$} for small \mbox{$t\neq 0$}.
\end{proof}

\begin{remark}\label{1.25}{\rm
The family $C_t$ of the above proof defines a deformation of the germ
$(C(q),z)$ 
to \mbox{$(C_t,\{z, z_t\})$}, where \mbox{$(C_t,z) \sim (C,z)$} and
\mbox{$(C_t,z_t)$} is an
ordinary $n$--fold point, $n$ as in Lemma \ref{1.20}.   In particular,
$(C(q), z)$ is a degeneration of a germ which is topologically equivalent to
$(C,z)$. 
}
\end{remark}

\section{$h^1$--vanishing criterion for
zero--dimensional schemes of class $\kg\ks_1$ in the plane}
\setcounter{equation}{0}

\begin{lemma}
 \label{2.1}
For any \mbox{$d\ge 1$} and \mbox{$X\in\kg\ks_1$} satisfying
\begin{equation}
  \label{2.2}
\deg X<(3-2\sqrt{2})(d-\mt \, X)^2
\end{equation}
there is \mbox{$X'\sim X$} with \mbox{$h^1(\kj_{X'/\P^2}(d))=0$}.
\end{lemma}

\begin{proof}
We shall prove the following statement.
Let \mbox{$L\subset\P^2$} be a fixed straight line. There exist
\mbox{$\alp,\bet\ge 0$} such that for any integer
\mbox{$d\ge 1$} and \mbox{$X\in\kg\ks_1$} satisfying
\begin{align}
\deg X & \,\le \; \bet(d-\mt \, X)^2 \label{2.3}\\
\deg(X\cap L) & \,\le \; d-\alp \,\frac{\deg X}{d} \label{2.4}
\end{align}
there exists \mbox{$X'\sim X$} with
\mbox{$\deg(X'\cap L)=\deg(X\cap L)$}, \mbox{$h^1(\kj_{X'/\P^2}(d))=0$}.
Moreover, in Step 2,  we show that for our approach the maximal possible value
for 
$\beta$ is attained at 
\begin{equation}
  \label{2.10}
\alp=\sqrt{2}+1,\quad \bet=3-2\sqrt{2}\: .
\end{equation}

{\it Step 1}. Assume that $X$ is an ordinary singularity, that
means \mbox{$T^\ast_X=\{z\}$}. Then the ideal of $X$ in
\mbox{$\ko_{\P^2,z}$} 
is defined by the vanishing of the coefficients of all monomials lying under
the diagonal
\mbox{$\,[(0,\mt \, X),(\mt \, X,0)]\,$} in the Newton diagram.
Since \mbox{$\mt \, X<d$} by (\ref{2.3}),
these (linear) conditions are independent, hence
\mbox{$h^1(\kj_{X/\P^2}(d))=0$}.

So, further on, we can suppose that \mbox{$\deg(X) >0$} and that $X$ is not
an ordinary singularity. We proceed by induction in $d$. For \mbox{$d\leq 2$}
there is nothing to consider. In the induction step, we reduce $X$ by $L$ and
have to show
\renewcommand{\arraystretch}{1.7}
\begin{equation}
  \label{2.7}
\begin{array}{c}
\deg(X:L)\;=\;\deg X-\deg(X\cap L)\;\le\;\bet(d-1-\mt \,(X:L))^2
\ ,\\
\deg((X:L)\cap L)\;\le\;\deg(X\cap L)\;\le\;
d-1-\alp\,\displaystyle{\frac{\deg(X:L)}{d-1}}
\: .
\end{array}
\end{equation}
Then, by 
the induction assumption
\mbox{$h^1(\kj_{Y/\P^2}(d-1))=0$} for some \mbox{$Y\sim X:L$},
\mbox{$\deg(Y\cap L)=\deg((X:L)\cap L)$}. By Proposition \ref{1.17a} there
exists \mbox{$X'\sim X$} with \mbox{$X':L=Y$} and \mbox{$\deg(X'\cap
  L)=\deg(X\cap L)$}.
Since \mbox{$h^1(\kj_{X\cap L/L}(d))=0$},
because \mbox{$\deg(X'\cap L)\le d-\alp\cdot\deg X/d<d+1$}, we obtain by
Proposition \ref{1.11} (iv)
the desired relation \mbox{$h^1(\kj_{X'/\P^2}(d))=0$}.

{\it Step 2}. Assume that
\begin{equation}
  \label{2.6}
\deg(X\cap L)\:=\:d-\alp \,\frac{\deg X}{d}\: .
\end{equation}

Due to \mbox{$\mt \,(X:L)\le\mt \, X$}, (\ref{2.3}) and (\ref{2.6}), the first
inequality in (\ref{2.7}) will follow from 
$$\bet(d-\mt \, X)^2-d+\alp \,\frac{\bet(d-\mt \, X)^2}{d}
\;\le\;\bet(d-1-\mt \, X)^2,$$
which is equivalent to
$$d^2(1-\alp\bet-2\bet)
+\alp\bet(2d-\mt \, X)\cdot\mt \, X+2\bet d\cdot\mt \, X+\bet d\;\ge\; 0,$$
hence, due to \mbox{$\mt \, X\le d$}, it is enough to impose the condition
\begin{equation}
  \label{2.8}
1\;\ge\;(\alp+2)\bet\: .
\end{equation}
The second inequality in
(\ref{2.7}) will follow from
$$(\alp+\alp^2)\cdot\deg X \;\le\; (\alp-1)d^2+d\: ,$$
which, by (\ref{2.3}), holds true as
\begin{equation}
  \label{2.9}
\alp-1\;\ge\;\bet(\alp+\alp^2).
\end{equation}

We are interested in $\bet$ as large as possible.
The inequality (\ref{2.9}) gives
$$\bet\;\le\; \frac{\alp-1}{\alp+\alp^2}\;\le\; 3-2\sqrt{2},$$
and the maximal value is attained at \mbox{$\alp=\sqrt{2}+1$}. So, from now on
we suppose (\ref{2.10}), especially the condition (\ref{2.8}) is satisfied.

{\it Step 3}. Assume that
$$\deg(X\cap L)\,<\,d-\alp\,\frac{\deg X}{d}\ ,$$
$X$ is not an ordinary singularity,
there exists a branch $Q$ of \mbox{$(C,z)$} such that
\mbox{$T^\ast_X\cap L=T^\ast_X\cap Q$}, and there is no branch $Q'$ of $(C,z)$
with
\mbox{$T^\ast_X\cap L\subsetneq T^\ast_X\cap Q'$}.

In this case \mbox{$T^\ast_X\cap L$} consists of at least two points. Therefore
\mbox{$\deg(X:L)<\deg X$} and
\mbox{$\deg((X:L)\cap L) \le \deg(X\cap L)-2$}. 
Moreover, by (\ref{2.3}) and
(\ref{2.8}), we have \mbox{$\,\alp \deg X /d\leq d\!-\!1$}, hence 
$$
\deg((X:L)\cap L)\; \le \;
(d-1)-\alp\,\displaystyle{\frac{\deg(X:L)}{d-1}}\:.
$$
On the other hand, by Lemma
\ref{1.15}, \mbox{$\mt \,(X:L)=\mt \, X-1$}, thus
$$\deg(X:L)\;< \;\deg X \;\le\; \bet(d-\mt \, X)^2
\;=\; \bet(d-1-\mt \,(X:L))^2.$$

{\it Step 4}. Assume that
$$\deg(X\cap L)<d-\alp\,\frac{\deg X}{d}\: ,$$
and there is a branch $Q$ of
\mbox{$(C,z)$} such that: (1) \mbox{$T^\ast_X \cap Q$} consists of points
\mbox{$z_1=z,\dots ,z_r$}, naturally ordered, (2)
\mbox{$T^\ast_X \cap L$} 
consists of points \mbox{$z_1,\dots ,z_s$}, \mbox{$1\leq s<r$}, (3) the
multiplicity
\mbox{$m=\mt \,(C_{(z_{s+1})},z_{s+1})$} satisfies
\begin{equation}
  \label{2.11}
\deg(X\cap L)+m\:>\:d-\alp\,\frac{\deg X}{d}\: .
\end{equation}
In this case, Lemma \ref{1.6} gives
\begin{equation}
  \label{2.12}
\deg X\,>\,\frac{(\mt \, X)^2+m^2}{2}\: .
\end{equation}
Since \mbox{$m\le\mt \, X$} and due to (\ref{2.11}), the first
inequality in (\ref{2.7}) will follow from
$$\bet(d-\mt \, X)^2-d+\alp\,\frac{\bet(d-\mt \, X)^2}{d}+\mt \, X \:\le\:
\bet(d-1-\mt \, X)^2,$$ 
or, equivalently, from
\begin{equation}
  \label{2.13}
(1-2\beta-\alpha\beta)\left(d-\mt \,X\right)+\beta  + \alpha
\beta \,\frac{\mt\, X}{d} \left(d- \mt \,X \right)\: \geq \:0\,,
\end{equation}
which holds true if (\ref{2.8}) is satisfied.
Similarly, the second inequality in
(\ref{2.7}) will follow from
$$d-\alp\,\frac{\deg X}{d} \;\le\; d-1- \alp\,\frac{\deg X-d+\alp\deg X/d+
m}{d-1},$$
which, by (\ref{2.3}), is satisfied, if
\begin{equation}
\label{2.15}
(\alp-1)d^2+d-\alp md \ge (\alp+\alp^2) \,\bet\,(d-\mt \, X)^2 \ .
\end{equation}
The coefficient of $d^2$ is zero by (\ref{2.10}), hence, it is enough to show
that 
$$2\bet(\alp\!+\!\alp^2)\,d\cdot\mt \, X-
\bet(\alp\!+\!\alp^2)(\mt \, X)^2-\alp md\;\ge\; 0\ ,$$
or, equivalently,
$$m\,\le\,\bet(1+\alp)\left(2\lam-\frac{\lam^2}{d}\right) \!\;=:\,
\varphi(\lam)\,, \qquad \lam=\mt \, X.$$
Indeed, due to (\ref{2.12}) an (\ref{2.3}), we have \mbox{$\lambda\leq
  \frac{\sqrt{2\beta}}{1+\sqrt{2\beta}}\,d$} and 
\begin{equation}
  \label{2.16}
m\,\le\,\psi(\lam)\,:=\,\left\{ 
\begin{array}{cl}
\lam &\mbox{for }\; 0\le\lam\le \frac{\sqrt{\bet}}{1+\sqrt{\bet}}\,d,\\
\sqrt{2\bet(d-\lam)^2-\lam^2} &\mbox{for }\;
\frac{\sqrt{\bet}}{1+\sqrt{\bet}}\,d
\le \lam\le \frac{\sqrt{2\bet}}{1+\sqrt{2\bet}}\,d\,.
\end{array}
\right.
\end{equation}
Since \mbox{$\varphi(d\sqrt{\bet}/(\sqrt{\bet}\!+\!1))=\psi(d\sqrt{\bet}/
(\sqrt{\bet}\!+\!1))$} as (\ref{2.10}) holds,
and $\varphi$ is increasing concavely in the segment
\mbox{$[0\,,\,d\sqrt{2\bet}/(\sqrt{2\bet}\!+\!1)]$}, we obtain
\mbox{$m\le\psi(\lam)\le \varphi(\lam)$}.

{\it Step 5}. Assume that
$$\deg(X\cap L)\:<\:d-\alp\,\frac{\deg X}{d}\ ,$$
and that there is a (smooth) branch $Q$ of
\mbox{$(C,z)$} such that: (1) \mbox{$T^\ast_X \cap Q$} consists of points
\mbox{$z_1=z,\dots,z_r$}, naturally ordered, (2) \mbox{$T^\ast_X \cap L$}
consists of points \mbox{$z_1,\dots,z_s$}, \mbox{$1\leq s<r$}, (3) the
multiplicity
\mbox{$m=\mt \,(C_{(z_{s+1})},z_{s+1})$} satisfies
$$\deg(X\cap L)+m\:\le\: d-\alp\,\frac{\deg X}{d}\ .$$
Then by Lemma \ref{1.14} we specialize the point $z_{s+1}$ on the line
$L$ and consider the new scheme \mbox{$\widetilde{X}\sim X$} with
\mbox{$\deg(\widetilde{X}\cap L)=\deg(X\cap L)+m$}. By the semi--continuity
of cohomology, \mbox{$h^1(\kj_{\widetilde{X}/\P^2}(d))=0$} yields
\mbox{$h^1(\kj_{X/\P^2}(d))=0$}. Thus, specializing points of $Q$ onto
$L$ we come to one of the cases studied above.
\end{proof}

\pagebreak[3]

\section{$H^1$--Vanishing Criterion for Zero--Dimensional Schemes of Class
$\kg\ks$ in the plane}
\setcounter{equation}{0}

For a scheme \mbox{$X\in\kg\ks$} denote by \mbox{$\mt_s X$}
the sum of the multiplicities of all singular branches
of the underlying germ \mbox{$(C,z)$}. Note that $\mt \, X$,
$\mt_s X$ are invariant with respect to the extension (cf.~(\ref{1.19})).

\begin{lemma}
  \label{3.1}
For any integer \mbox{$d\ge 1$} and any \mbox{$X\in\kg\ks$} satisfying
\begin{equation}
  \label{3.2}
\deg X\:\le\:\bet_0(d-\mt \, X-\mt_s X)^2\ ,
\end{equation}
where \mbox{$\bet_0=(\alp_0+8)^{-1}=0.10340\dots$} and
\mbox{$\alp_0=(31-3\sqrt{85})/2 =1.6706\dots$}
is the positive root of the equation
$$\left(\frac{\sqrt{4\alp^3+\alp^2-4\alp}+\alp-2}{2(1+\alp+\alp^2)}
\right)^2=\frac{1}{\alp+8}\ ,$$
there is \mbox{$X'\sim X$} with
\mbox{$h^1(\kj_{X'/\P^2}(d))=0$}.
\end{lemma}

\begin{proof}
As in the proof of Lemma \ref{2.1}, we shall obtain a more general statement.
Let \mbox{$L\subset\P^2$} be a fixed straight line. There exist
\mbox{$\alp,\bet>0$} such that, for any integer \mbox{$d\ge 1$} and any
\mbox{$X\in\kg\ks$}, satisfying
$$\deg X\:\le\:\bet(d-\mt \, X-\mt_s X)^2,\quad
\deg(X\cap L)\:\le \:d-\alp\,\frac{\deg X}{d}\ ,$$
there exists \mbox{$X'\sim X$} with
\mbox{$h^1(\kj_{X'/\P^2}(d))=0$}, \mbox{$\deg(X'\cap L)=\deg(X\cap L)$}.
Finally we show that for $\alp$, $\bet$ we can take the values $\alpha_0$ and
$\bet_0$, respectively. 

{\it Step 1}. In the case \mbox{$X\in\kg\ks_1$} the proof of Lemma \ref{2.1}
gives sufficient conditions on $\alp$, $\bet$, namely (\ref{2.8}), (\ref{2.9})
in the Steps 2, 3, and the inequality (\ref{2.15}) in Step 4. Due to
(\ref{2.16}), it is sufficient to check the inequality (\ref{2.15})
after removing the term $d$ and substituting
\mbox{$d\sqrt{\bet}/(\sqrt{\bet}+1)$} for $\mt \, X$ and $m$, or, equivalently,
to have
\begin{equation}
  \label{3.3}
\bet\:\le\:\left(\frac{\sqrt{4\alp^3+\alp^2-4\alp}+\alp-2}{
2(1+\alp+\alp^2)}\right)^2\ .
\end{equation}

{\it Step 2}. Assume that \mbox{$X\in\kg\ks \backslash \kg\ks_1$}.
Since \mbox{$\mt_s X\le\mt \, X$}, we can perform the inductive procedure
described in the proof of Lemma \ref{2.1} under assumptions (\ref{2.8}),
(\ref{2.9}), (\ref{3.3}), until the following situation occurs:
\begin{itemize}
\item[(1)] $X$ satisfies (\ref{3.2});
\item[(2)] let \mbox{$T^\ast_X\cap L=\{q_1,\dots ,q_N\}$},
  such that  for each branch $D$ going through \mbox{$q:=q_N$} 
  $$ \mt \,(D_{(q_1)},q_1)=\ldots = \mt \,(D_{(q_{N-1})},q_{N-1})>
  \mt \,(D_{(q)},q)>0$$
(especially $D$ is not a smooth branch), and 
\begin{equation}
  \label{3.4}
\deg(X\cap L)\:<\:d-\alp\,\frac{\deg X}{d}-m_q\ ,
\end{equation}
where $m_q$ is the multiplicity of $X$ at $q$.
\end{itemize}

Again, in the induction step, it is sufficient to show the two inequalities
\renewcommand{\arraystretch}{1.7}
\begin{equation}
  \label{3.4A}
\begin{array}{c}
\deg X-\deg(X\cap L)\;\le\;\bet(d-1-\mt \,(X:L)-\mt_s(X:L))^2
\ ,\\
\deg((X:L)\cap L)\;\le\;\deg(X\cap L)\;\le\;
d-1-\alp\,\displaystyle{\frac{\deg(X:L)}{d-1}}
\: .
\end{array}
\end{equation}

Consider the possible situations.

{\it Step 3}. Under the hypotheses of the second step, assume that
$$\deg(X\cap L)\;\ge\; d-2m'-\alp\,\frac{\deg X}{d-m'}\ ,$$
where $m'$ is the sum of the multiplicities of all branches, going
through $q$.

Then the first inequality in (\ref{3.4A}) will follow from
$$\alp\,\frac{\bet(d\!-\!\mt \, X\!-\!\mt_s X)^2}{d\!-\!m'}\:\le\:
(d\!-\!2m')+\bet-2\bet(d\!-\!\mt \, X\!-\!\mt_sX)\ .$$
Since \mbox{$m'\le\mt_s X \le \mt \, X$}, replacing the left--hand side by
\mbox{$\alp\bet(d-\mt \, X-\mt_sX)$}, and replacing the term
\mbox{$d-2m'$} in the right--hand side by \mbox{$d-\mt \, X-\mt_sX$}, one
obtains a stronger inequality, namely
$$0\:\le\:\bet+(1-2\bet-\alp\bet)(d-\mt \, X-\mt_sX)\ ,$$
which is an immediate consequence of (\ref{2.8}).

The second inequality of (\ref{3.4A}) will follow from
$$d-1-\alp\,\frac{\deg X-\deg(X\cap L)}{d-1}\;\ge\;\deg(X\cap L)\ .$$
We replace $\deg X$ and $\deg(X\cap L)$ by the upper bounds (\ref{3.2}),
(\ref{3.4}) and obtain
$$(\alp-1)d^2+d\;\ge\;\bet(\alp+\alp^2)(d-\mt \, X-\mt_sX)^2-
m_q(d-1-\alp)\ ,$$
which holds true by (\ref{2.9}).

{\it Step 4}. Under the hypotheses of the second step, assume that
\begin{equation}
  \label{3.5}
\deg(X\cap L)\:<\:d-2m'-\alp\,\frac{\deg X}{d-m'}\ .
\end{equation}
In this case we have to exert ourselves to obtain an analogue to the first
inequality in (\ref{3.4A}). For that, we shall perform the following $m'$--step
algorithm.

Let \mbox{$1\le j\le m'$} and let
\mbox{$X_{j-1}\in\kg\ks$} be defined by a germ \mbox{$(C_{j-1},z)$} and
a tree $T^\ast_{j-1}$ (\mbox{$X_0=X,\ C_0=C,\ T^*_0=T^*$}), such that
at the endpoint $q$ of \mbox{$T^\ast_{j-1}\cap L$} the strict transform of
$C_{j-1}$ has the multiplicity
$m_q^{(j-1)}$ (\mbox{$m_q^{(0)}=m_q$}), and
$$\deg(X_{j-1}\cap L)\:<\:d-m'-\alp\,\frac{\deg X}{d\!-\!m'}\ .$$
The $j$-th step of the algorithm appears as follows: introduce
$$s_j:=\min\left\{l\ge 0\ \Big|\ \deg(X_{j-1}\cap L)+lm'_j\,\ge\,
d-2m'-\alp\,\frac{\deg X}{d\!-\!m'}\right\},$$
where $m'_j$ is the sum of the multiplicities of all branches of
\mbox{$(C_{j-1},z)$} going through $q$ at the preceding point
\mbox{$\overline{q} \in T^\ast_{j-1}\cap L$}. In particular, 
\mbox{$m'_1=m'$} and \mbox{$s_1\geq 1$}. Define $X'_{j-1}$ as the extension
$$X'_{j-1}:=X_{j-1}\underbrace{(q)\dots(q)}_{s_j\ \mbox{\scriptsize times}}\ ,
\:\mbox{ and }\: X_j:=X'_{j-1}:L\ .$$

Note that in the previous formula, in the definition of $X_{j-1}$ and
in the assumption of Step 2, we denote different points by $q$. But all
these points appear in the extension operation introduced above, and the
notation $q$ moves to new points of new schemes as was described
in the assertion of Lemma \ref{1.20}.

Due to Lemma \ref{1.22}, again it is enough to show 
\begin{align}
\deg(X_{m'}\cap L) &\: \le \; d-m'-\alp\,\frac{\deg X_{m'}}{d-m'}\
,\label{3.8}\\ 
\deg X_{m'} & \:\le\; \bet(d-m'-\mt \, X_{m'}-\mt_sX_{m'})^2 \label{3.9}
\end{align}
to complete the induction step.



\smallskip\noindent
We define the set 
$$\Lambda\,:=\:  \bigl\{\!\:j\in[1,m'\!-\!\!\:1]\:\big|\: m'_{j+1}<m'_{j}
\!\:\bigr\} 
\:=\: \bigl\{\!\:j_1,j_2,\dots, j_\ell\!\:\bigr\}\,,$$
\mbox{$j_{k+1}> j_{k}$}, that is, by Proposition \ref{1.11} (ii), we have
$$ m'_{j_k+1}\: =\: m'_{j_k+2}\:=\ldots =\:
m'_{j_{k+1}}\:=\:m'\!-\!\!\:k 
$$
for any \mbox{$k=0,\dots,\ell$} (where \mbox{$j_0:=0$} and
\mbox{$j_{\ell+1}:=m'$}). We set   
$$ N_k \,:=\, \sum_{i=j_{k}+1}^{j_{k+1}} s_{i}\,. $$ 
Note that if \mbox{$j\in \Lambda$}, that is, if \mbox{$m'_{j+1}=m'_{j}-1$}, 
then there are two possibilities:  first, it might be that \mbox{$\mt\, X_j=\mt
\,X_{j-1}\!\!\:-\!\!\:1\ $}.  Secondly, if this is not the case, then 
Lemma \ref{1.16} (ii) gives at least the existence of a branch $Q'$ of
the germ \mbox{$(C_j,z)$} with \mbox{$T^\ast_j\cap Q'\subset T^\ast_j\cap L$},
and Lemma \ref{1.15} implies that \mbox{$\mt \,X_{j+1}=\mt\,
X_{j}\!\!\:-\!\!\:1$}. 
In any case, we have \mbox{$\mt_s X_j=\mt_s X_{j-1}\!\!\:-\!\!\:1\ .$}
Hence, we can estimate
$$ \widetilde{\ell}\::=\: (\mt\, X\!\!\:+\!\!\:\mt_s X) - (\mt\,
X_{m'}\!\!\:+\!\!\: \mt_s X_{m'}) \:\left\{
\renewcommand{\arraystretch}{1.1}
\begin{array}{ll} 
\geq 0 & \text { if $\ell = 0$}\,,\\
\geq \ell\!\!\;+\!\!\;1 & \text { if $\ell\neq 0$} \,,
\end{array}
\right. $$
To run an induction step, it is sufficient to
show that
\begin{equation}
\label{3.14}
\deg X_{m'} \,\leq \, \deg X\,, \qquad 
\deg X_{m'} \, \leq \, \beta
(d\!\!\:-\!\!\:m'\!+\!\!\:\widetilde{\ell}\!\!\:-\!\!\:\mt\, 
X\!\!\:-\!\!\:\mt_s X)^2.
\end{equation}
By construction, we have
\begin{eqnarray*}
\deg X_{m'} & = & \deg X + \sum_{k=0}^{\ell} N_k
\!\;\frac{(m'\!-\!\!\:k)(m'\!-\!\!\:k\!\!\:+\!\!\:1)}{2} \!\;- 
\sum_{j=0}^{m'-1}  
\deg(X'_{j}\cap L) \\
& = & \deg X +
\frac{(m'\!-\!\!\:\ell\!\;)(m'\!-\!\!\:\ell\!\!\;+\!\!\:1)}{2} 
\cdot \sum_{k=0}^{\ell} N_k  \!\:- \sum_{j=0}^{m'-1}  \deg(X'_{j}\cap L) \\
&& \phantom{\deg X } +
\sum_{k=0}^{\ell-1} \bigl((m'\!-\!\!\:k)\cdot (N_0+\ldots + N_{k})\bigr) \,,
\end{eqnarray*}
and we can estimate
\begin{equation}
\label{partial sum}
 \deg(X_{j_{k+1}-1}\cap L) \: \geq \: (m' \!-\!\!\:k) \cdot 
(1+N_0+\ldots + N_{k}) \,,
\end{equation}
for any \mbox{$k=0,\dots,\ell\!\!\;-\!\!\:1$}. On the other hand, in the $j$-th
step we have 
$$d\!\!\;-\!\!\;2m'\!+\!\!\;m'_{j}\!\!\;-\!\!\;\alpha\,\frac{\deg
X}{d\!\!\:-\!\!\:m'} 
\:>\: \deg(X'_{j-1}\cap L) \: \geq \: 
d\!\!\;-\!\!\;2m'\!-\!\!\;\alpha\,\frac{\deg X}{d\!\!\:-\!\!\:m'} \,,$$
\mbox{$j=0,\dots,m'\!-\!\!\:1$}. In particular, this together with
(\ref{partial sum}) implies that
$$ (m' \!-\!\!\:\ell\!\:) \cdot \sum_{k=0}^{\ell} N_k \:\leq \: 
 \deg(X_{m'-1}\cap L) -(m'\!-\!\!\: \ell) \: < \:
 d\!\!\;-\!\!\;2m'\!-\!\!\;\alpha\,\frac{\deg X}{d\!\!\:-\!\!\:m'} \,.
$$
Hence, we can estimate 
\begin{eqnarray*}
\deg X_{m'} & \leq  & \deg X +
\frac{m'\!-\!\!\:\ell\!\!\;+\!\!\:1}{2}  
\cdot (m'\!-\!\!\:\ell\!\;)\cdot \sum_{k=0}^{\ell} N_k  \!\:- 
\sum_{\renewcommand{\arraystretch}{0.5}
\begin{array}{c}
\scriptstyle{ k=0 } \\
\scriptstyle{ k+1\not \in \Lambda} 
\end{array}
}^{m'-1}  \deg(X'_{k}\cap L)\\
& \leq &  \deg X + \left(\frac{m'\!-\!\!\:\ell\!\!\;+\!\!\:1}{2} -
m'\!\!\:+\!\!\;\ell\!\:\right) \cdot  
\left(d\!\!\;-\!\!\;2m'\!-\!\!\;\alpha\,\frac{\deg X}{d\!\!\:-\!\!\:m'}
\right)\\
& \leq & \deg X - \frac{m'\!-\!\!\:1\!\!\:-\!\!\;\ell}{2} \cdot
\left(d\!\!\;-\!\!\;2m'\!-\!\!\;\alpha\,\frac{\deg X}{d\!\!\:-\!\!\:m'} 
\right) \,,
\end{eqnarray*}
whence the first inequality in (\ref{3.14}). The second will
follow from
\begin{eqnarray*}
\lefteqn{\beta(d\!\!\:-\!\!\:m'\!+\!\!\:\widetilde{\ell}\!\!\:-\!\!\:\mt\,
X\!\!\:-\!\!\:\mt_s X)^2}\hspace{2cm}\\ 
&\geq & \beta(d\!\!\:-\!\!\:\mt\, X\!-\!\!\:\mt_s X)^2-
\frac{m'\!-\!\!\:1\!\!\:-\!\!\;\ell}{2} \cdot
\left(d\!\!\;-\!\!\;2m'\!-\!\!\;\alpha\,\frac{\deg X}{d\!\!\:-\!\!\:m'} 
\right)
\end{eqnarray*}
or, equivalently,
\begin{eqnarray*}
\lefteqn{
2\beta(m'\!-\!\!\;\widetilde{\ell}\,)
(d\!\!\:-\!\!\:\mt\, X\!-\!\!\:\mt_s X)}\hspace{2cm}\\
 &\leq &
\frac{m'\!-\!\!\:1\!-\ell}{2}\cdot \left(d-2m'\!\!\:-
\alpha\beta\,\frac{(d\!\!\:-\!\!\:\mt\, X\!-\!\!\:\mt_s X)^2}{
d\!\!\:-\!\!\:m'}\right) +(m'\!-\!\!\;\widetilde{\ell}\,)^2. 
\end{eqnarray*}
Since \mbox{$m'\le\mt_s X \le \mt\, X $}, this holds whenever
\begin{equation}
  \label{3.23}
0 \: \leq \: \frac{1\!\!\:-\!\!\:
\alpha\beta}{4\beta} \cdot (m'\!-\!\!\:1\!\!\:-\!\!\;\ell)(d\!\!\:-\!\!\:2m')-
(m'\!-\!\!\;\widetilde{\ell}\,)(d\!\!\:-\!\!\:2m')\, .
\end{equation}
Note that for fixed $\ell$ the right-hand side takes its 
minimum for the minimal possible value of $\widetilde{\ell}$. Hence, it
suffices to consider two cases:

\smallskip \noindent
{\it Case 4A.} \mbox{$\,\widetilde{\ell}=\ell=0\!\;$}. Then, since
\mbox{$m'\!\geq 2$},  (\ref{3.23}) is
implied by 
\begin{equation}
\label{3.13}
\frac{1}{\beta} \: \geq \: 8+\alpha\,.
\end{equation}

\smallskip \noindent
{\it Case 4B.} \mbox{$\,\widetilde{\ell}=\ell+1\geq 2\!\;$}. Then (\ref{3.23})
holds if
\mbox{$\frac{1}{\beta}  \geq  4+\alpha$},
which is a consequence of (\ref{3.13}).


{\it Step 5}. Finally we look for the maximal value of $\bet$,
satisfying (\ref{2.8}), (\ref{2.9}), (\ref{3.3}), (\ref{3.13}) with any
\mbox{$\alp>1$}.
It is easily shown that $\bet_0$ mentioned in the assertion of
Lemma \ref{3.1} is this maximal value.

{\it Step 6}. The induction base for \mbox{$d\le 6$} is trivial. Indeed, for
\mbox{$X\in \kg\ks\backslash \kg\ks_1$}
the right--hand side of (\ref{3.2}) does not exceed
\mbox{$\bet_0(6-2-2)^2<1$}. 

\end{proof}

\section{Existence of Curves with one Singular Point}
\setcounter{equation}{0}

We start with the following auxiliary statement:

\begin{lemma}
  \label{4.1}
Let a scheme \mbox{$Y\in\kg\ks\cap\ks$} be defined by a
germ \mbox{$(C,z)$} with branches \mbox{$Q_1,\dots,Q_p$} and the tree
\mbox{$T^\ast=T^\ast(C)$}.
Let \mbox{$f\in J_Y$} satisfy
\begin{equation}
  \label{5.2}
\mt \,(f,z)\,=\,\mt \, Y,\quad(f, Q_i)\,>\!\!\sum_{q\in T^*\!\cap
  Q_i}\!m_q\cdot \mt \,(Q_{i,(q)},q),\ \,i=1,\dots,p.
\end{equation}
Then the germ \mbox{$(f,z)$} also defines the singularity scheme $Y$,
in particular
\mbox{$(f,z)$} and \mbox{$(C,z)$} have the same topological type (which we
denote by $Y$).
\end{lemma}

\begin{proof}
{\it Step 1}. In the case \mbox{$Y\in\kg\ks_1\cap\ks$} this
easily can be shown by induction on $\deg Y$.
The induction base with \mbox{$Y=\emptyset$} and $C$ being non--singular at
$z$, is trivial. Hence, assume \mbox{$\deg Y>0$} and blow up the point $z$.
If $w$ is an intersection point of the strict transform $C^\ast$ of $C$
with the exceptional divisor $E$, then, for any branch \mbox{$Q_{i,(w)}$} of
$C^\ast$ centred at $w$, we have
\begin{equation}
  \label{4.2}
(f_{(w)},Q_{i,(w)})\:=\:(f\cdot Q_i)-\mt \,(f,z) \:>\!
\sum_{q\in T^\ast\cap Q_i\backslash \{z\} }\!\!\!m_q\:\ge\: 0\, .
\end{equation}
Thus, \mbox{$w\in f_{(w)}$} and, if $C^\ast$ is non--singular at $w$, then
\mbox{$\mt \,(f_{(w)},w)\ge 1=\mt \,(C^\ast\!,w)$}. If $C^\ast$ is singular at
$w$, then by Lemma \ref{1.8}, \mbox{$f_{(w)}\in J_{Y^\ast_{(w)}}$}, where
\mbox{$Y^\ast_{(w)}\in\kg\ks_1$}  is defined by the germ
\mbox{$(C^\ast\!,w)$} and its tree of essential points. Especially, again,  
\mbox{$\mt \,(f_{(w)},w)\ge\mt \,(C^\ast\!,w)$}. On the other hand,
$$\mt \,(f,z)\:\ge\:\sum_{w\in C^\ast\cap E}\!\!\mt
  \,(f_{(w)},w)\:\ge\!\sum_{w\in C^\ast\cap 
  E}\! \mt \,(C^\ast\!,w)\:=\:(C^\ast\!,E)\:=\:\mt \, Y,$$
which implies \mbox{$\mt \,(f_{(w)},w)=\mt \,(C^\ast\!,w)$} and, together with
(\ref{4.2}) and the induction assumption,
  \mbox{$(f_{(w)},w)$} defines the same singularity scheme as
\mbox{$(C^\ast\!,w)$}. Since, moreover, $C^\ast$ and $f_{(w)}$ are
transversal to $E$ at $w$, \mbox{$(f,z)$} and \mbox{$(C,z)$} define the same
singularity scheme $Y$.
 
{\it Remark}. Let \mbox{$(C,z)$} be given in local coordinates $x$,$y$ at $z$
by 
$$ C(x,y)\:=\: \prod_{i=1}^m \left(y-\sum_{j=1}^\infty a_{ij}x^j\right) $$
and define the {\bf essential part} $\Gam_{es}$ of the Newton diagram $\Gam$ of
\mbox{$C(x,y)$} at the origin as the union of
\begin{itemize}
\item[(i)] all the integral points \mbox{$(i,j)\in \Gam$} with positive
  $i,\,j$, 
\item[(ii)] a point \mbox{$(n,0)\in\Gam$}, if it is not an endpoint of an edge
\mbox{$[(n,0),(n',1)]\subset\Gam$},
\item[(iii)] a point \mbox{$(0,n)\in\Gam$}, if it is not an endpoint of an edge
\mbox{$[(0,n),(1,n')]\subset\Gam$}.
\end{itemize}
We claim that \mbox{$f(x,y)$} has the same
essential part $\Gam_{es}$ of the Newton diagram at the origin.
This is easily shown by induction on $\deg Y$,
using the transformation
\begin{equation}
  \label{4.3}
(x,y)\mapsto(x,xy)
\end{equation}
as the blowing--up at $z$.

{\it Step 2}.
Now assume that $Y$ is arbitrary in $\kg\ks\cap \ks$.
We apply induction on the number
$$N\,:=\,\sum_Q\sum_{q\in T^\ast\!\cap Q}\!\mt \,(Q_{(q)},q),$$
where $Q$ runs through all singular branches of $C$. The case
\mbox{$N=0$} means just \mbox{$Y\in\kg\ks_1$}.
If \mbox{$N>0$} then, again, we blow up the point $z$. 
Each intersection point of the strict transform $C^\ast$ of $C$
with the exceptional divisor $E$ corresponds to a straight line $W$
through $z$, tangent to $C$. Without restriction, we can suppose that in local
coordinates $x,y$ at $z$ we have $W=y$ and $C(x,y)$ decomposes in local
branches
\mbox{$Q_1,\dots,Q_p$} with
$$ Q_k=\prod_{i=1}^{s_k}
\left( y-\xi_i^{(k)} (x)\right),\;\;\;\xi_i^{(k)} (x)=\sum_{j=1}^\infty
a^{(k)}_{ij}x^{j/s_k},\;\;\; s_k=\mt \,(Q_k,z), $$  
The (covering) transformation 
\begin{equation}
  \label{4.4}
(x,y)\mapsto(x^M\!,y) \;\mbox{ with } \,M:=\prod_{k=1}^p s_k
\end{equation}
turns \mbox{$(C,z)$} into a germ \mbox{$(\widetilde{C},z)$} with
multiplicity \mbox{$\mt \,(\widetilde{C},z)=\mt \,(C,z)=\mt \, Y$}
and only non--singular branches
$$Q_i^{(k)}\:=\: y-\sum_{j=1}^{\infty}a^{(k)}_{ij}x^{jM/s_k},\quad
i=1,\dots,s_k,\; k=1,\dots,p\,.$$
Let \mbox{$\widetilde{Y}\in\kg\ks_1$} be defined by the germ
\mbox{$(\widetilde{C},z)$} and the tree \mbox{$T^\ast(\widetilde{C})$}.
We shall show that the transform \mbox{$\widetilde{\varphi}(x,y):=
\varphi(x^M\!,y)$} of any element \mbox{$\varphi\in J_Y$} belongs to
\mbox{$J_{\widetilde{Y}}$}.

By Lemma \ref{1.8}, we have for any (fixed) \mbox{$i=1,\dots,s_k$},
$$
\left(\widetilde{\varphi}, Q_i^{(k)}\right) 
 \: = \: \displaystyle\frac{M}{s_k}(\varphi, Q_k)
\:\ge\: \displaystyle\frac{M}{s_k}\left(2\del(Q_k)+\sum_{l\ne k}(Q_l,
 Q_k)+\!\!\sum_{q\in T^\ast\!\cap Q_k}\!\!\!\mt \,(Q_{k,(q)},q)\right)\, .
$$
Hence, using the considerations in the proof of Lemma \ref{1.17},
\begin{align*}
\left(\widetilde{\varphi}, Q_i^{(k)}\right) 
&\: \geq \; \left(\sum_{j\neq i} (Q_i^{(k)},Q_j^{(k)})-M+\frac{M}{s_i}\right) +
\left( \sum_{l\neq k} \sum_{j=1}^{s_l} (Q_i^{(k)},Q_j^{(l)}) \right)\\
&\quad\quad\quad\quad + \displaystyle \left( \max_{(j,l)\neq (i,k)}
(Q_i^{(k)}, Q_j^{(l)})+ M - \frac{M}{s_i} \right)\\
&\:=\;  \displaystyle \sum_{(j,l)\ne(i,k)} (Q_i^{(k)},
Q_j^{(l)})\,+\,\# \,\{ q\in T^\ast(\widetilde{C})\cap Q_i^{(k)} \} \, ,
\end{align*}
and Lemma \ref{1.8} implies \mbox{$\widetilde{\varphi}\in J_{\widetilde{Y}}$}.

The germ \mbox{$\widetilde{f}(x,y)=f(x^M\!,y)$}, clearly, satisfies
\mbox{$\mt \,(\widetilde{f},z)=\mt \, Y$}. Due to (\ref{5.2}), the previous
computation with \mbox{$\widetilde{f}$} instead of $\widetilde{\varphi}$
gives
$$(\widetilde{f},Q_i^{(k)})\;>\!\!\sum_{q\in T^\ast\!(\widetilde{C})\cap
Q_i^{(k)}}\!\!\mt \,(\widetilde{C}_{(q)},q),\quad i=1,\dots,s_k,\;
k=1,\dots,p\,.$$
Hence, by Step 1, \mbox{$(\widetilde{f},z)$} defines the same singularity
scheme as \mbox{$(\widetilde{C},z)$}.
Denote by $\Gam$ the Newton diagram of \mbox{$C(x,y)$} at the
origin. Evidently, the Newton diagram $\widetilde{\Gam}$ of $\widetilde{C}$ and
its essential part $\widetilde{\Gam}_{es}$ are
obtained from
$\Gam$, $\Gam_{es}$ by the transformation \mbox{$(I,J)\mapsto(MI,J)$}. As
established above, $\widetilde{f}$ has the same essential part
$\widetilde{\Gam}_{es}$ of the Newton diagram at the origin.
Therefore, $\Gam_{es}$ is the essential part of the Newton diagram
of $f(x,y)$ at the origin.

Let $\Gam'$ be the part of $\Gam$ corresponding to
the branches of $C$ tangent to $y$, and
$$(n_1,m_1),\dots,(n_l,m_l),\quad m_1>\dots>m_l=0,$$
be the vertices of $\Gam'$. Applying the blowing--up (\ref{4.3}) at $z$,
we easily obtain that the Newton diagram
of \mbox{$C^\ast(x,y)$} at \mbox{$w=(0,0)$} has the vertices
$$(0,m_1),\ (n_2\!-\!n_1,m_2),\dots,\ (n_{l-1}\!-\!n_1,m_{l-1}),\
(n_l\!-\!n_1,0),$$
and that \mbox{$f_{(w)}(x,y)$} has the Newton diagram with vertices
$$(0,m_1),\ (n_2\!-\!n_1,m_2),\dots,\ (n_{l-1}\!-\!n_1,m_{l-1}),\ (r,0),$$
where $r$ may be different from \mbox{$n_l\!-\!n_1$} only in the case
\mbox{$m_{l-1}=1$}. In particular, this means that
\begin{equation}
  \label{4.9}
\mt \,(f_{(w)},w)=\mt \,(C^\ast\!,w),\quad(f_{(w)}, E)=\mt \,(C^\ast\!\cap
E,w)\, .
\end{equation}
On the other hand, for any branch \mbox{$Q_{i,(w)}$} of $C^\ast$
centred at $w$, we have
$$(f_{(w)},Q_{i,(w)})\,=\,(f, Q_i)-\mt \,(f,z)\cdot\mt \, (Q_i,z)\,>\!\!
\sum_{q\in T^\ast\!\cap Q_i\backslash \{z\}}\!\!m_q\cdot\mt \,(Q_{i,(q)},q)\
.$$
This, together with (\ref{4.9}) and the induction assumption, implies that
\mbox{$(f_{(w)}\cdot E,w)$} defines the same singularity scheme as
\mbox{$(C^\ast\!\cdot E,w)$}.

Finally, blowing down all the germs $f_{(w)}$, \mbox{$w\in C^\ast\!\cap E$},
one obtains that \mbox{$(f,z)$} also defines the singularity scheme $Y$.
\end{proof}

\begin{definition}
  \label{4.10}
Let $F$ be a curve of degree $d$ with an isolated singular
point $z$. The germ \mbox{$(H_{\P^2}^{es,d},F)$} of the equisingular stratum
\mbox{$H_{\P^2}^{es,d} \subset \P^N=\P(\Gam(\ko_{\P^2}(d)))$},
\mbox{$N=d(d\!+\!3)/2$}, at $F$
(cf.~\cite{GrL}) is called {\bf T--smooth}, if it is smooth and, for any
\mbox{$d'>d$}, it is a transversal intersection in
\mbox{$\P(\Gam(\ko_{\P^2}(d')))$} of \mbox{$(H_{\P^2}^{es,d'}\!,F)$}
and
\mbox{$\P(\Gam(\ko_{\P^2}(d)))$}, included in
\mbox{$\P(\Gam(\ko_{\P^2}(d')))$} by \mbox{$C\mapsto CL^{d'-d}$}, where $L$ is
a fixed generic straight line not passing through $z$. 
\end{definition}

\begin{lemma}
  \label{4.11}
For any scheme \mbox{$X\in\kg\ks\cap\ks$},
and any positive integer $d$, satisfying
\begin{equation}
  \label{4.12}
\deg X+\mt \, X+1\:<\:\left\{
\begin{array}{cl}
(3-2\sqrt{2})(d-\mt \, X-2)^2\,,&\mbox{if}
\ X\in\kg\ks_1,\\
\bet_0(d-\mt \, X-\mt_sX-2)^2\,,&\mbox{if}
\ X\in\kg\ks\backslash\kg\ks_1,
\end{array}
\right.
\end{equation}
$\bet_0$ as above in Lemma \ref{3.1},
there exists an irreducible curve $F$ of degree $d$ with a singular
point $z$ of (topological) type $X$ as its only singularity such that
the germ \mbox{$(H_{\P^2}^{es,d},F)$} is T--smooth. 
\end{lemma}

\begin{proof}
{\it Step 1}. Let \mbox{$(C,z)$} be a defining germ of the scheme $X$ and
consider the germ \mbox{$(\widetilde{C},z)$}, \mbox{$\widetilde{C}=C\cdot L$},
where $L$ is a generic straight line through $z$.
Note that
$$T^\ast(\widetilde{C})=T^\ast(C),\quad\mt \,(\widetilde{C},z)=\mt \,(C,z)+1\
.$$
Now we introduce (1) the scheme $X'$, defined by the
germ \mbox{$(C,z)$} and the tree $T^\ast_{X'}$, containing $T^\ast(C)$ and
the first non--essential points of all local branches of \mbox{$(C,z)$}, and
(2) the scheme $\widetilde{X}$, defined by the germ \mbox{$(\widetilde{C},z)$}
  and the tree \mbox{$T^\ast(C)$}.
Clearly,
\begin{equation}
  \label{4.13}
  \begin{split}
&\mt \, X'=\,\mt \, X,\quad\mt \,\widetilde{X}=\,\mt \, X+1,\\
\deg X'\leq&\,\deg X+\mt \, X,\quad\deg\widetilde{X}\leq\,
\deg X+\mt \, X+1,\\
  \end{split}
\end{equation}
and in addition, for any local branch $Q$ of \mbox{$(C,z)$}
and any elements \mbox{$f\in\kj_{X'/\P^2,z}$},
\mbox{$g\in\kj_{\widetilde X/\P^2,z}$},
\begin{equation}
  \label{4.14}
  \begin{split}
&(f, Q)\:>\!\sum_{q\in T^\ast(C)\cap Q}\mt \,(C_{(q)},q),\\
(g, Q)\:\ge&\sum_{q\in T^\ast\!(C)\cap Q}\!\!\mt \,(\widetilde{C}_{(q)},q)\:>\!
\sum_{q\in T^\ast\!(C)\cap Q}\!\!\mt \,(C_{(q)},q).\\
  \end{split}
\end{equation}

By (\ref{4.12}), (\ref{4.13}) and the Lemmas \ref{2.1},
\ref{3.1}, we may assume that
\begin{equation}
\label{4.16a}
h^1(\kj_{X'/\P^2}(d-1))=h^1(\kj_{\widetilde{X}/\P^2}(d-1))=0.
\end{equation}
Hence, due to \mbox{$X'\subsetneq\widetilde{X}$}, there exists a (generic)
curve 
\begin{equation}
  \label{4.15}
f\in H^0(\kj_{X'/\P^2}(d\!-\!1))\backslash H^0(\kj_{\widetilde{X}
/\P^2}(d\!-\!1)),
\end{equation}
which, by (\ref{4.13}), (\ref{4.14}), satisfies the condition of Lemma
\ref{4.1}. Thus, \mbox{$(f,z)$} defines the singularity scheme $X$.
Replacing, if
necessary, multiple components of $f$ (which do not go through $z$)
by distinct components, we obtain a reduced curve $f$ of degree \mbox{$d-1$}.

If $f$ is irreducible, then we are done. Otherwise, we shall use Bertini's
Theorem to
construct the desired irreducible curve $F$ of degree $d$. For that, let a
straight line $L$ meet $f$ at \mbox{$d-1$} distinct
non--singular points $w_1,\dots,w_{d-1}$. Obviously,
\mbox{$h^1(\kj_{\{w_i\}/L}(1))=0$} for each \mbox{$i=1,\dots,d\!-\!1$} and, by
(\ref{4.16a}) \mbox{$h^1(\kj_{X'/f}(d\!-\!1))=0$}. First, observe that this
implies
\begin{equation}
\label{4.18}
h^1(\kj_{X'\cup\{w_i\}/fL}(d))=0,\quad 1\le i\le d\!-\!1.
\end{equation}
Indeed, the first morphism \mbox{$\mbox{id}_1\otimes L+f \otimes
\mbox{id}_2$} in the exact sequence
$$0\longrightarrow\ko_{f}(d-1)\oplus\ko_{L}(1)\longrightarrow\ko_{fL}(d)
\longrightarrow \ko_{f\cap L}\longrightarrow 0$$
maps the sheaf \mbox{$\kj_{X'/f}(d-1)\oplus\kj_{\{w_i\}/L}(1)$}
injectively to the sheaf \mbox{$\kj_{X'\cup\{w_i\}/fL}(d)$}. Now, consider
the commutative diagram
$$
\arraycolsep0.1cm
\begin{array}{ccccccccc}
\ &\ &\ &\ & 0 &\ &\ &\ &\ \\
\ &\ &\ &\ & \downarrow &\ &\ &\ &\ \\
\ &\ &\ &\ &\ko_{\P^2} &\ &\ &\ &\ \\
\ &\ &\ &\ & \downarrow &\ &\ &\ &\ \\
0 &\longrightarrow &\kj_{X'\cup\{w_i\}/\P^2}(d)
&\longrightarrow &\ko_{\P^2}(d) &\longrightarrow
&\ko_{X'\cup\{w_i\}} &\longrightarrow & 0 \\
\ &\ &\ &\ & \downarrow &\ &\| &\ &\ \\
0 &\longrightarrow &\kj_{X'\cup\{w_i\}/fL}(d)
&\longrightarrow &\ko_{fL}(d)
&\longrightarrow
&\ko_{X'\cup\{w_i\}} &\longrightarrow & 0 \\
\ &\ &\ &\ & \downarrow &\ &\ &\ &\ \\
\ &\ &\ &\ & 0 &\ &\ &\ &\ 
\end{array}
$$
to deduce from (\ref{4.18}) the surjectivity of
\mbox{$H^0(\ko_{\P^2}(d))\to H^0(\ko_{X'\cup\{w_i\}})$}, which is equivalent to
$$
h^1(\kj_{X'\cup\{w_i\}/\P^2}(d))=0,\quad 1\le i\le d\!-\!1.
$$
In particular, there exist curves
$$\varphi_i\in H^0(\kj_{X'/\P^2}(d))\backslash H^0(\kj_{X'\cup\{w_i\}/
\P^2}(d)),
\quad 1\le i\le d\!-\!1,$$
and, by Bertini's theorem, the generic member $\Phi$ of the linear family
$$f\cdot L+\lam_1\varphi_1+\dots+\lam_{d-1}\varphi_{d-1}$$
is irreducible and smooth outside of the base points, whereas at $z$ it has
a singularity of the (topological) type $X$. 
Assume that $\Phi$ has singular points \mbox{$w_1,\dots,w_p$} different
from $z$. Using the preceding
arguments, we show that there exist curves
$$\Phi_j\in H^0(\kj_{X'/\P^2}(d))\backslash H^0(\kj_{X'\cup\{w_j\}/\P^2} (d)),
\quad j=1,\dots,p.$$
Finally, by Bertini's theorem, a generic member $F$ of the linear family
$$f\cdot L+\lam_1\varphi_1+\dots+\lam_{d-1}\varphi_{d-1}
+\mu_1\Phi_1+\dots+\mu_p\Phi_p$$
has no singularities outside $z$.

{\it Step 2}. Note that, by Lemma \ref{1.7}, \mbox{$h^1(\kj_{X'/\P^2}(d))=0$}
implies  \mbox{$h^1(\kj^{\mbox{\scriptsize es}}(d))=0$}, where
\mbox{$\kj^{\mbox{\scriptsize es}}\subset\ko_{\P^2}$} is the ideal sheaf of the
zero--dimensional scheme given by the equisingularity ideal
\mbox{$I^{\mbox{\scriptsize 
      es}}\subset \ko_{\P^2,z}$} of \mbox{$(F,z)$}. But the latter equality
yields the required T--smoothness of the germ \mbox{$(H_{\P^2}^{es,d},F)$}
(\cite{GrK,GrL}). 
\end{proof}

\section{Proof of Theorems 1 and 2}
\setcounter{equation}{0}

\begin{definition}
  \label{5.1}
Let \mbox{$X\in \ks$} be a singularity scheme. By $s(X)$ we denote the minimal
integer 
such that there exists an irreducible curve $F$ of degree $s(X)$
with a singular point $z$ of type $X$ as its only singularity, and such that
the germ \mbox{$(H_{\P^2}^{es,d},F)$} of the equisingular stratum
\mbox{$H_{\P^2}^{es,d} \subset \P^N=\P(\Gam(\ko_{\P^2}(d)))$},
\mbox{$N=d(d\!+\!3)/2$}, at $F$ is T--smooth.
\end{definition}

\begin{lemma}
\label{5.1a} 
Let \mbox{$X\in \ks$} be a singularity scheme. Then $s(X)\leq \sig(X)$, where
$$\sig(X):=
\left\{
\renewcommand{\arraystretch}{1.7}
\begin{array}{cl}
\left[\sqrt{2} \,\mt \, X\right] +1, &\mbox{ if }\ X \mbox{ is ordinary,}\\
\displaystyle \left[(1\!+\!\sqrt{2})\sqrt{\deg X\!+\!\mt \,
    X\!+\!1}\right]\!+\!\mt \, X\!+\!3,
&\mbox{ if } X\in\kg\ks_1,\\
\left[\sqrt{(\deg X\!+\!\mt \, X\!+\!1)\bet_0^{-1}}\right]\!+\!\mt \,
X\!+\!\mt_sX\!+\!3, & \mbox{ if }
X\not\in\kg\ks_1.
\end{array}
\right.
$$
\end{lemma}

\begin{proof}
By Lemma \ref{4.11}, \mbox{$s(X)\le\sig(X)$} for any
non--ordinary singularity $X$. The case of an ordinary singularity was already
treated in \cite{GLS} and the result follows from Lemma \ref{5.3} below.
\end{proof}

\begin{definition}
  \label{5.2a}
Let \mbox{$\oz=(z_1,\dots,z_n)$} be a tuple of distinct points in $\P^2$,
and \mbox{$\om=(m_1,\dots,m_n)$} be a vector of positive integers, then we
denote by
\mbox{$X(\oz,\om)$} the zero--dimensional scheme in $\P^2$ defined by
the ideals \mbox{$\mathfrak{m}_{z_i}^{m_i}\subset\ko_{\P^2,z_i}$},
\mbox{$i=1,\dots,n$}. 
\end{definition}

\begin{lemma}[GLS, section 3.3]
  \label{5.3}
Let \mbox{$z_1,\dots,z_n$} be distinct generic points in $\P^2$, and the
positive integers \mbox{$d,m_1,\dots,m_n$} satisfy
$$\sum_{i=1}^n\frac{m_i(m_i+1)}{2}\:<\:\frac{d^2+6d-1}{4}-
\left[\frac{d}{2}\right]\, .$$
Then (1) there exists an irreducible curve $F_d$ of degree $d$ with ordinary
singular points \mbox{$z_1,\dots,z_n$} of multiplicities
\mbox{$m_1,\dots,m_n$},
respectively, as its only singularities; (2) for \mbox{$\oz=(z_1,\dots,z_n)$},
 \mbox{$\om=(m_1,\dots,m_n)$},
\begin{equation}
  \label{5.4}
h^1(\kj_{X(\oz,\om)/ \P^2}(d))=0\, .
\end{equation}
\end{lemma}

\begin{lemma}
  \label{5.5}
If a positive integer $d$ and the singularities \mbox{$S_1,\dots,S_n$}
satisfy the inequality
$$\sum_{i=1}^n \frac{(s(S_i)+1)(s(S_i)+2)}{2}\:<\: \frac{d^2+6d-1}{4}-\left[
\frac{d}{2}\right],$$
then there exists an irreducible curve $F$ of degree $d$ with
exactly $n$ singular points of (topological) types \mbox{$S_1,\dots,S_n$},
respectively.
\end{lemma}

\begin{proof}
Due to Lefschetz's principle (cf.~\cite{JL}, Theorem 1.13), we suppose, without
loss of generality, \mbox{$\K=\C$}.
Let
$$G_i(x,y)=\sum_{0\le k+l\le s(S_i)}a^{(i)}_{kl}
x^ky^l,\quad \deg G_i=s(S_i),\quad i=1,\dots,n,$$
be affine irreducible curves such that each of
them has exactly one singular point at \mbox{$(0,0)$} of type
\mbox{$S_1,\dots,S_n$}, respectively.

On the other hand, let \mbox{$z_1,\dots,z_n$} be distinct generic points in
$\P^2$, then Lemma \ref{5.3} implies the existence of an irreducible curve
$G'$ of degree $d$ having ordinary singularities at \mbox{$z_1,\dots,z_n$} of
multiplicities \mbox{$m_i=s(S_i)+1$}, \mbox{$i=1,\dots,n$}, as its only
singularities.
For any \mbox{$i=1,\dots,n$}, let us fix affine coordinates $x_i$, $y_i$ in a
neighbourhood of $z_i$.
The relation (\ref{5.4}) means
that an affine neighbourhood $U$ of $G'$ in
\mbox{$\P^N=\P(\Gam(\ko_{\P^2}(d)))$}
can be parametrized
by the following independent parameters: to any \mbox{$\Phi\in U$}
we assign
\begin{itemize}
\itemsep0.1cm
\item[(1)] coefficients \mbox{$A^{(i)}_{kl}$}, \mbox{$0\le k+l\le s(S_i)$},
  from the representation
$$\Phi(x_i,y_i)=\sum_{0\le k+l\le d}\! A^{(i)}_{kl}x^k_iy^l_i$$
of $\Phi$ in the coordinates $x_i$, $y_i$, for any \mbox{$i=1,\dots,n$},
\item[(2)] some parameters $B_j$, 
\mbox{$1\le j\le r:=N-\sum_{i=1}^n (s(S_i)\!+\!1)(s(S_i)\!+\!2)/2$}.
\end{itemize}
First we deform $G'$ into a curve $G$ by addition of the leading form, that is,
the part of degree $s(S_i)$, of
\mbox{$G_i(x,y)$} as the $s(S_i)$--jet at $z_i$, \mbox{$i=1,\dots,n$}. Without
loss of generality we may suppose that $G$ is irreducible and has no
singularities outside \mbox{$\{z_1,\dots,z_n\}$}. Moreover, in the local
coordinates  $x_i$, $y_i$, $G$ is represented as
$$G^{(i)}(x_i,y_i)=\sum_{d\ge k+l\ge s(S_i)}\!a^{(i)}_{kl}x_i^ky_i^l,$$
\mbox{$i=1,\dots,n$}, and corresponds to the parameter values
\mbox{$A^{(i)}_{kl} = 0$}, if \mbox{$k+l<s(S_i)$}, \mbox{$A^{(i)}_{kl}
= a^{(i)}_{kl}$}, if \mbox{$k+l=s(S_i)$}, and \mbox{$B_j=b_j$},
\mbox{$j=1,\dots,r$}.
We shall look for the desired curve $F$ close to $G$, given by parameters
$$
A^{(i)}_{kl}(\tau)  =  \left\{ 
\begin{array}{cr}
\tau^{s(S_i)-k-l}a^{(i)}_{kl}(\tau) & \mbox{ if } \,k+l<s(S_i)\\
a^{(i)}_{kl}(\tau) & \mbox{ if } \,k+l=s(S_i)
\end{array}
\right.
\; \mbox{ and } \; B_j=b_j(\tau)\,,
$$
where \mbox{$a^{(i)}_{kl}(\tau)$}, \mbox{$b_j(\tau)$}
are smooth functions in a neighbourhood of zero such that
$$a^{(i)}_{kl}(0)\,=\,a^{(i)}_{kl},\quad b_j(0)\,=\,b_j,$$
for all \mbox{$i,j,k,l$}. 

Let us fix some \mbox{$1\leq i\leq n$}. In a neighbourhood of the point
$z_i$ we have
$$F(x_i,y_i)=\sum_{0\le k+l\le s(S_i)}\!\!\tau^{s(S_i)-k-l}
a^{(i)}_{kl}(\tau)x_i^ky_i^l+\!\!
\sum_{d\geq p+q>s(S_i)}\!a^{(i)}_{pq}(\tau)x_i^py_i^q\,,$$
where the coefficients \mbox{$a^{(i)}_{pq}(\tau)$}, \mbox{$p+q>s(S_i)$}, are
affine functions in the parameters \mbox{$A^{(s)}_{kl}$}, \mbox{$s\neq i$},
and $B_j$, \mbox{$j=1,\dots,r$}.
The transformation \mbox{$(x_i,y_i)\mapsto (\tau x_i,\tau y_i)$} 
turns $F$ for sufficiently small \mbox{$\tau \neq 0$} into a curve
$$F_i(x_i,y_i)=\sum_{0\le k+l\le s(S_i)}\! a^{(i)}_{kl}(\tau)
x_i^ky^l_i\,+\!\!\sum_{d\geq p+q>s(S_i)} \!\tau^{p+q-s(S_i)}a^{(i)}_{pq} (\tau)
x_i^py_i^q$$
close to $G_i$ in \mbox{$\P^N=\P(\Gam(\ko_{\P^2}(d)))$}.
By the definition of $s(S_i)$,
the germ of the equisingular stratum \mbox{$H_{\P^2}^{es,d}\subset \P^N$} at
$G_i$ can be described by \mbox{$c(S_i)$} equations
$$\varphi^{(i)}_u(A^{(i)}_{kl})=0,\quad u=1,\dots,c(S_i),$$
on the coefficients \mbox{$A^{(i)}_{kl}$}, \mbox{$0\le k+l\le d$}, of a curve
$$H(x_i,y_i)=\sum_{0\le k+l\le d} \!A^{(i)}_{kl}x_i^ky_i^l,$$
such that there exists
\mbox{$\Lam_i\subset\{(k,l)\in\Z^2\,|\: k,l\ge 0, \ k+l\le s(S_i)\}$}, \mbox{$
\mbox{card}(\Lam_i)=c(S_i)$},
with
$$\mbox{det}\left(\frac{\partial\varphi^{(i)}_u}{\partial
A^{(i)}_{kl}}\right)_{
\renewcommand{\arraystretch}{0.6}
\begin{array}{c}
\scriptstyle{u=1,\dots,c(S_i)}\\
\scriptstyle{(k,l)\in\Lam_i}
\end{array}}
\ne 0$$
at the point
\mbox{$A^{(i)}_{kl}=a^{(i)}_{kl}$}, \mbox{$0\le k\!+\!l\le s(S_i)$} and 
\mbox{$ A^{(i)}_{pq}=0$}, \mbox{$ p\!+\!q>s(S_i)$}.
Thus, the condition on $F$ to have singular points of types
\mbox{$S_1,\dots,S_n$} can be expressed as the system of equations
\begin{equation}
  \label{5.6}
\varphi^{(i)}_u \left( \{a^{(i)}_{kl}(\tau)\,|\:k\!+\!l\le s(S_i)\},\; \{
  \tau^{p+q-s(S_i)}a^{(i)}_{pq}(\tau)\,|\: p\!+\!q > s(S_i)\} \right) =0\,, 
\end{equation}
\mbox{$u=1,\dots,c(S_i)$}, \mbox{$i=1,\dots,n$}.
From the above it follows immediately that at the point
\mbox{$\tau=0$} the determinant
$$\mbox{det} \left( \frac{\partial \varphi^{(i)}_u}{
\partial a^{(m)}_{kl}}\right)_{
\renewcommand{\arraystretch}{0.6}
\begin{array}{c}
\scriptstyle{u=1,\dots,c(S_i)}\\
\scriptstyle{(k,l)\in\Lam_m}\\
\scriptstyle{m,i=1,\dots,n}
\end{array}}
=\prod_{i=1}^n
\mbox{det}\left(\frac{\partial
\varphi^{(i)}_u}{\partial a^{(i)}_{kl}}\right)_{
\renewcommand{\arraystretch}{0.6}
\begin{array}{c}
\scriptstyle{u=1,\dots,c(S_i)}\\
\scriptstyle{(k,l)\in\Lam_i}
\end{array}}
$$
does not vanish, which implies the existence of an appropriate solution
to (\ref{5.6}), and, hereby, the existence of a curve $F$ with
$n$ singular points of types \mbox{$S_1,\dots,S_n$}.

The only thing we should explain, is why $F$ has no other singular points.
Let us consider $F$ as a polynomial function.
First, note that $F$ is a small deformation of the function $G$,
which has \mbox{$(d\!-\!1)^2-(s(S_1)\!-\!1)^2-\dots -(s(S_n)\!-\!1)^2$}
critical points out of the zero level. Second, we deform each
ordinary critical point $z_i$ of $G$ by means of the
function $G_i$ which has \mbox{$(s(S_i)\!-\!1)^2-\mu(S_i)$} critical
points out of the zero level. Hence, $F$ has at least
\mbox{$(d\!-\!1)^2-\mu(S_1)-...-\mu(S_n)$} critical points out of
the zero level, that means the $n$  constructed singular points
are the only singular points of the curve $F$.
\end{proof}

Thus, to prove the Theorems 1 and 2 from the introduction, by Lemma \ref{5.1a},
it suffices to prove the following

\begin{proposition}
  \label{5.8}
For any singularity $X\in \ks$,
$$(\sig(X)\!+\!1)(\sig(X)\!+\!2)\le 196\,\mu(X)\, .$$
\end{proposition}

\begin{proof}
If $X$ is an ordinary singularity, then \mbox{$\mu(X)=(\mt \, X\!-\!1)^2$};
therefore 
$$(\sig(X)\!+\!1)(\sig(X)\!+\!2)\,\leq \,2(\mt \, X\!+\!1)(\mt \,
X\!+\!3)\,\leq\, 30\,\mu(X).$$ 
For an arbitrary singularity $X\in \ks$, defined by the germ \mbox{$(C,z)$}, we
have 
$$\mu(X)\,=\,2\del(X)-r+1\,=\,\sum_{q\in T^\ast\!(C)} m_q(m_q\!-\!1)-r+1,$$
where $r$ is the number of local branches.
Note that the number of the essential points $q$ with \mbox{$m_q=1$} does not
exceed \mbox{$\mt \, X<\sqrt{\mu(X)}+1$}. Hence, by Lemma \ref{1.6},
$$
2(\deg X+\mt \, X+1) \: \leq \: \sum_{m_q>1} m_q(m_q\!+\!1) +4\mt \, X+2 \:< \:
3\left(\sqrt{\mu(X)}+\sqrt{2}\right)^2\,; 
$$
thus
$$
\sig(X) \: < \: \sqrt{ \frac{\deg X+\mt \, X+1}{\bet_0}}+2\mt \, X+3
\: < \:
\left(
  \sqrt{\frac{3}{2\bet_0}}+2\right)\sqrt{\mu(X)}+\sqrt{\frac{3}{\bet_0}}+5\,, 
$$
and, finally, it follows that \mbox{$(\sig(X)\!+\!1)(\sig(X)\!+\!2)$} is
smaller than
$$
\left(\sqrt{\frac{3}{2\bet_0}}
\!+\!2\right)^2\mu(X)+
\left(\sqrt{\frac{3}{\bet_0}}\!+\!6\right)\left(\sqrt{\frac{3}{\bet_0}}
\!+\!7\right)
+\left(2 \sqrt{\frac{3}{\bet_0}}\!+\!13 \right)
\left(\sqrt{\frac{3}{2\bet_0}}\!+\!2\right)\sqrt{\mu(X)},
$$
which for \mbox{$\mu(X)\ge 2$} does not exceed \mbox{$196\mu(X)$}.
\end{proof}


\end{document}